\pdfoutput=1
\newif\ifpersonal
\newif\ifarxiv
\personaltrue % comment to remove personal notes
\arxivtrue % comment to display shortened version for journal submissions
\RequirePackage[l2tabu, orthodox]{nag} %detect whether obsolete packages are used
\documentclass[10pt,a4paper,reqno,oneside]{amsart} %reqno places equation numbers on the right
\usepackage{calligra}
\usepackage{amsmath,amsthm,amssymb,mathrsfs,mathtools,bm,eucal,tensor} % math related
\usepackage{microtype,fixltx2e} % latex technical issues
\usepackage[scaled]{beramono,berasans}
\usepackage{enumerate,comment,braket,xspace,tikz-cd, etaremune} %utilities
\usepackage[all,cmtip]{xy} % utilities
\usepackage[utf8]{inputenc} % input encoding
\usepackage[T1]{fontenc} % font encoding
\usepackage{lmodern}
\definecolor{linkcolor}{HTML}{005050}
\usepackage[centering,vscale=0.7,hscale=0.7]{geometry}
\usepackage{hyperref}
\usepackage[capitalize]{cleveref}
\usepackage{capt-of}
\usepackage{graphicx}
\usepackage{xparse}
\usepackage{url}
\usepackage[toc,page]{appendix}

\usepackage{vmargin}
\setpapersize{A4}
\setmarginsrb{25mm}{10mm}{25mm}{10mm}%
{12mm}{10mm}{5mm}{10mm}

\usepackage{fancyhdr}
\theoremstyle{plain}
\newtheorem{thm-intro}{Theorem}
\newtheorem{thm}{Theorem}[section]

\newtheorem*{thm*}{Theorem}
\newtheorem{lem}[thm]{Lemma}

\newtheorem*{lem*}{Lemma}
\newtheorem{prop}[thm]{Proposition}

\newtheorem{Ansatz}[thm]{Ansatz}
\newtheorem{Hypothesis}[thm]{Hypothesis}
\newtheorem{proposition}[thm]{Proposition}

\newtheorem{conj}[thm]{Conjecture}
\newtheorem{Question}[thm]{Question}
\newtheorem{cor}[thm]{Corollary}
\newtheorem{cor-intro}[thm-intro]{Corollary}

\theoremstyle{definition}
\newtheorem{definition}[thm]{Definition}
\newtheorem{defin}[thm]{Definition}

\newtheorem{defin-intro}[thm-intro]{Definition}

\theoremstyle{remark}
\newtheorem{rem}{Remark}

\newtheorem{example}[thm]{Example}

\newtheorem{eg-intro}[thm-intro]{Example}
\newtheorem{rem-intro}[thm-intro]{Remark}
\numberwithin{equation}{section}

\usepackage[backend=biber, style=numeric, sorting = nyt , giveninits=true,doi=false,url=false,maxbibnames=99]{biblatex}
\ifpersonal
\newcommand*{\personal}[1]{\textcolor[rgb]{0.6,0.6,1}{(Personal: #1)}}
\newcommand*{\todo}[1]{\textcolor{red}{(Todo: #1)}}
\else
\newcommand*{\personal}[1]{\ignorespaces}
\newcommand*{\todo}[1]{\ignorespaces}
\fi

\newcommand{\C}{\mathbb C}

\newcommand{\Q}{\mathbb Q}
\newcommand{\R}{\mathbb R}
\newcommand{\Z}{\mathbb Z}
\newcommand{\N}{\mathbb N}

\newcommand{\bbR}{\mathbb R}

\newcommand{\sw}{\mathsf w}

\newcommand{\fX}{\mathfrak X}

\newcommand{\cK}{\mathcal K}

\newcommand{\cO}{\mathcal O}

\DeclareFontFamily{U}{BOONDOX-calo}{\skewchar\font=45 }
\DeclareFontShape{U}{BOONDOX-calo}{m}{n}{<-> s*[1.05] BOONDOX-r-calo}{}
\DeclareFontShape{U}{BOONDOX-calo}{b}{n}{<-> s*[1.05] BOONDOX-b-calo}{}
\DeclareMathAlphabet{\mathcalboondox}{U}{BOONDOX-calo}{m}{n}
\newcommand{\bbA}{\mathbb A}

\newcommand{\bbP}{\mathbb P}

\newcommand{\bbH}{\mathbb H}

\newcommand{\Gr}{\text{Gr}}

\makeatletter
\let\save@mathaccent\mathaccent
\newcommand*\if@single[3]{%
	\setbox0\hbox{${\mathaccent"0362{#1}}^H$}%
	\setbox2\hbox{${\mathaccent"0362{\kern0pt#1}}^H$}%
	\ifdim\ht0=\ht2 #3\else #2\fi
}
\newcommand*\rel@kern[1]{\kern#1\dimexpr\macc@kerna}
\newcommand*\widebar[1]{\@ifnextchar^{{\wide@bar{#1}{0}}}{\wide@bar{#1}{1}}}
\newcommand*\wide@bar[2]{\if@single{#1}{\wide@bar@{#1}{#2}{1}}{\wide@bar@{#1}{#2}{2}}}
\newcommand*\wide@bar@[3]{%
	\begingroup
	\def\mathaccent##1##2{%
		\let\mathaccent\save@mathaccent
		\if#32 \let\macc@nucleus\first@char \fi
		\setbox\z@\hbox{$\macc@style{\macc@nucleus}_{}$}%
		\setbox\tw@\hbox{$\macc@style{\macc@nucleus}{}_{}$}%
		\dimen@\wd\tw@
		\advance\dimen@-\wd\z@
		\divide\dimen@ 3
		\@tempdima\wd\tw@
		\advance\@tempdima-\scriptspace
		\divide\@tempdima 10
		\advance\dimen@-\@tempdima
		\ifdim\dimen@>\z@ \dimen@0pt\fi
		\rel@kern{0.6}\kern-\dimen@
		\if#31
		\overline{\rel@kern{-0.6}\kern\dimen@\macc@nucleus\rel@kern{0.4}\kern\dimen@}%
		\advance\dimen@0.4\dimexpr\macc@kerna
		\let\final@kern#2%
		\ifdim\dimen@<\z@ \let\final@kern1\fi
		\if\final@kern1 \kern-\dimen@\fi
		\else
		\overline{\rel@kern{-0.6}\kern\dimen@#1}%
		\fi
	}%
	\macc@depth\@ne
	\let\math@bgroup\@empty \let\math@egroup\macc@set@skewchar
	\mathsurround\z@ \frozen@everymath{\mathgroup\macc@group\relax}%
	\macc@set@skewchar\relax
	\let\mathaccentV\macc@nested@a
	\if#31
	\macc@nested@a\relax111{#1}%
	\else
	\def\gobble@till@marker##1\endmarker{}%
	\futurelet\first@char\gobble@till@marker#1\endmarker
	\ifcat\noexpand\first@char A\else
	\def\first@char{}%
	\fi
	\macc@nested@a\relax111{\first@char}%
	\fi
	\endgroup
}
\makeatother

\newcommand{\Coh}{\mathrm{Coh}}

\usetikzlibrary{decorations.markings} %arrows for open immersions and closed immersions
\tikzset{
  closed/.style = {decoration = {markings, mark = at position 0.5 with { \node[transform shape, xscale = .8, yscale=.4] {/}; } }, postaction = {decorate} },
  open/.style = {decoration = {markings, mark = at position 0.5 with { \node[transform shape, scale = .7] {$\circ$}; } }, postaction = {decorate} }
}

\DeclareMathOperator{\Aut}{Aut}

\DeclareMathOperator{\Bl}{Bl}

\DeclareMathOperator{\Hom}{Hom}
\DeclareMathOperator{\Image}{Im}

\DeclareMathOperator{\Pic}{Pic}

\DeclareMathOperator{\Spec}{Spec}

\DeclareMathOperator{\supp}{supp}

\DeclareMathOperator{\Symp}{Symp}

\title{Semi-stable degenerations of Calabi-Yau manifolds and mirror P=W conjecture}

\author{Sukjoo Lee}

\address{Department of Mathematics, University of Edinburgh, EH9 3FD, UK}
\email{Sukjoo.Lee@ed.ac.uk}

\begin{document}

\maketitle
\begin{abstract}
    % We propose a topological mirror construction of a Calabi-Yau manifold that admits a semi-stable degeneration of general type. This extends the work of Doran-Harder-Thompson where the type II case (Tyurin degeneration) was first investigated. The main ingredient is the notion of hybrid Landau-Ginzburg model, a multi-potential analogue of classical Landau-Ginzburg model, that serves as a mirror of each irreducible component of the degeneration fiber. This construction also provides a topological fibration to the projective space. 

    Mirror symmetry for a semi-stable degeneration of a Calabi-Yau manifold was first investigated by Doran-Harder-Thompson when the degeneration fiber is a union of two (quasi)-Fano manifolds. They propose a topological construction of a mirror Calabi-Yau that is a gluing of two Landau-Ginzburg models mirror to those Fano manifolds. We extend this construction to a general type semi-stable degeneration. As each component in the degeneration fiber comes with the simple normal crossing anti-canonical divisor, one needs the notion of a hybrid Landau-Ginzburg model, a multi-potential analogue of classical Landau-Ginzburg models. We show that these hybrid LG models can be glued to provide a topological mirror candidate of the Calabi-Yau which is also equipped with the fibration over $\bbP^N$. Furthermore, it is predicted that the perverse Leray filtration associated to this fibration is mirror to the monodromy weight filtration on the degeneration side \cite{MirrorClemensSchmid}. We explain how this can be deduced from the original mirror P=W conjecture \cite{PWmirror}.
    
\end{abstract}
\tableofcontents

\section{Introduction}
Traditionally, mirror symmetry is a conjectural relationship between two compact Kähler $n$-dimensional Calabi-Yau manifolds $X$ and $X^\vee$: the complex (algebraic) geometry of $X$ (B-side) is equivalent to the symplectic geometry of $X^\vee$ (A-side) and vice versa \cite{MirrorBook}\cite{MirrorBook2Clay}. The geometric construction of such pair was proposed by Strominger-Yau-Zaslow \cite{SYZ} by viewing them as a dual special Lagrangian torus fibration. This idea leads to extend mirror symmetry  beyond the Calabi-Yau cases, especially to the case of quasi-Fano manifolds \cite{AurouxFanoSYZ}. We say $X$ is quasi-Fano if the anti-canonical linear system $|-K_X|$ is base-point free and $H^i(X, \cO_X)=0$ for $i>0$. In this case, a mirror is given by so-called \textit{Landau-Gizburg (LG) model} $(Y, \omega, \sw:Y \to \C)$ where $(Y,\omega)$ is a $n$-dimensional Calabi-Yau Kähler manifold and $\sw:Y \to \C$ is locally trivial symplectic fibration near the infinity. We refer the readers to \cite{AurouxFanoSYZ} for more details.

\subsection{Generalization of Doran-Harder-Thompson construction}
A relationship between two different kinds of mirror symmetries for Calabi-Yau manifolds and quasi-Fano manifolds was first addressed by Doran-Harder-Thompson in the case of Tyruin degenerations \cite{mirrorsymmetrytyurin}. Recall that Tyurin degeneration is a semi-stable degeneration of a Calabi-Yau manifold $X$ into a union of two quasi-Fano varieties $X_1 \cup X_2$ over the smooth anti-canonical hypersurface $X_{12}:=X_1 \cap X_2$. This degeneration restricts the behavior of normal bundles of $X_{12}$ in $X_1$ and $X_2$ to be inverse to each other. Suppose one has a mirror LG model $(Y_i, \sw_i)$ for each pair $(X_i, X_{12})$ and generic fibers of each $\sw_i$ are topologically the same. It is predicted that the anti-canonical divisor $-K_{X_i}$ is mirror to the monodromy of a generic fiber $\sw_i^{-1}(t)$ near the infinity. By the adjunction formula, the normal bundle condition corresponds to the condition that the monodomies of $\sw_i^{-1}(t)$ are inverse to each other. This allows one to topologically glue two LG models to obtain a mirror candidate of $X$ which is also equipped with the map to $\bbP^1$. One natural question is how to generalize this construction when $X$ degenerates into a more general simple normal crossing variety. 

\begin{Question}
    How to extend the construction of Doran-Harder-Thompson for a semi-stable degeneration of more general types?
\end{Question}

The aim of this article is to answer this question and study related topics. It is predicted that for a Calabi-Yau mirror pair $X$ and $Y$, there is a correspondence between semi-stable degeneration of $X$ and the morphism from $Y$ to a projective variety (See Section \ref{sec:Toric degeneration}). To see this, we use recently developed language of hybrid Landau-Ginzburg models \cite{mypapermirror} that is a multi-potential analogue of classical LG model, whose idea goes back to \cite[Section 5.3]{AurouxFanoSYZLag}. We consider a triple $(Y, \omega, h:Y \to \C^N)$, called a \textit{hybrid Landau-Ginzburg model}, where $(Y,\omega)$ is a Calabi-Yau kähler manifold of dimension $n$ and $h=(h_1, \dots, h_N):Y \to \C^N$ is a fibration which is locally trivial around the infinity boundaries of the base (See Definition \ref{def: weak hybrid LG model general}). In fact, this turns out to be a suitable model to capture mirror symmetry of the quasi-Fano pair $(X,D)$ where $D=\cup_{i=1}^N D_i$ has $N$ irreducible components in the following sense: for any $I \subset \{1, \dots, N\}$, the induced quasi-Fano pair $(D_I:=\cap_{i \in I}D_i, \cup_{j \notin I}D_j\cap D_I)$ is expected to be mirror to the induced hybrid LG model $(Y_I:=\cap_{i \in I} \sw_i^{-1}(t_i), \omega|_{Y_I}, h|_{Y_I}:Y_I \to \C^{N-|I|})$ where $\sw_i^{-1}(t_i)$ is a generic fiber of $\sw_i$. We will review the precise notion of the hybrid LG model and the associated mirror symmetry relations in Section \ref{sec:extended Fano}.

Let's consider a semistable degeneration of a Calabi-Yau manifold $X$ into a simple normal crossing variety $X_c=\cup_{i=0}^NX_i$ whose dual boundary complex is the standard $N$-simples. Suppose we have a hybrid LG model $(Y_i, \omega_i, h_i:Y_i \to \C^N)$ mirror to each pair $(X_i, \cup_{j \neq i}X_{ij})$ with additional topological conditions (Hypothesis \ref{Hypothesis}). Similar to the Tyurin degeneration case, the semi-stability corresponds to the condition on the monodromies associated to the hybrid LG models (See Ansatz \ref{Ansatz:linebundle monodromy}). By shrinking the base of $h_i$ to the polydisk $\Delta_{h_i}$, this allows to topologically glue hybrid LG models and produces a symplectic fibration $\pi:Y \to \bbP^N$ (Proposition \ref{prop:Gluing LG models}) where the base $\Delta_{h_i}$ is identified with the locus $\{|z_j| \leq |z_i| |j=0, \dots, N\} \subset \bbP^N$. We also take a general hyperplane $H \subset \bbP^N$ and its complement $\bbP^N\setminus H \cong \C^N$. We write the induced fibration $\tilde{\pi}:\tilde{Y} \to \C^N$ for $\tilde{Y}:=\pi^{-1}(\bbP^N\setminus H)$. 

\begin{thm}\label{thm: topological mirror intro}(Theorem \ref{thm: topological mirror symmetry})
    Suppose that $(X_i, \cup_{j \neq i}X_{ij})$ is topological mirror to $(Y_i, h_i:Y_i \to \C^N)$ for all $i$. Then
    \begin{enumerate}
    \item $Y$ is topological mirror to $X$. In other words, $e(Y)=(-1)^ne(X)$
    \item $\tilde{Y}$ is topological mirror to $X_c$. In other words, $e(\tilde{Y})=(-1)^ne(X)$
    \end{enumerate}
    where $e(-)$ is the Euler characteristic. 
\end{thm}

\subsection{Mirror P=W conjecture}

The topological mirror relation in Theorem \ref{thm: topological mirror intro} is the weakest form of the mirror symmetry one would expect. We could ask further about other kinds of mirror symmetry relations for this construction. In the degeneration picture, there is a geometric auto-equivalence on the cohomology of $X$ induced by the monodromy of $X$ around the center $X_c$. It gives rise to the \textit{monodromy weight filtration} on $H^*(X)$ that constitutes the limiting mixed Hodge structure. On the other hand, on the cohomology of the central fiber $X_c$, there is Deligne's canonical weight filtration that constitutes the mixed Hodge structure. The natural question is what the corresponding filtrations on the mirror $Y$ and $\tilde{Y}$ are. 
\begin{Question}
    What is the filtration on the cohomology of $Y$ (resp. $\tilde{Y}$) that is mirror to the monodromy weight filtration (resp. the Deligne's canonical weight filtration)?
\end{Question}
The answer is expected to be the perverse Leray filtration associated to $\pi:Y \to \bbP^N$ (resp. $\tilde{\pi}:\tilde{Y} \to \C^N$), proposed by Doran-Thompson \cite[Conjecture 4.3]{MirrorClemensSchmid}. 
\begin{conj}\label{conj: mirror PW degeneration intro}
 \begin{enumerate}
     \item For $X$ and $Y$, we have 
     \begin{equation}
        \dim_\C(\Gr_F^p\Gr^{W_{\lim}}_{p+q}H^{p+q+l}(X,\C))=\dim_\C(\Gr_F^{n-q}\Gr^{P}_{n+p-q}H^{n+p-q+l}(Y,\C))
    \end{equation}
    where $P_\bullet$ is the perverse Leray filtration associated to $\pi$.
    \item For $X_c$ and $\tilde{Y}$, we have
    \begin{equation}
        \dim_\C(\Gr_F^p \Gr^W_{p+q} H^{p+q+l}(X_c,\C))=\dim_\C(\Gr_F^{n-q}\Gr^{P}_{n+p-q} H_c^{n+p-q+l}(\tilde{Y}),\C)
    \end{equation}
    where $P_\bullet$ is the perverse Leray filtration associated to $\tilde{\pi}$.
 \end{enumerate}
 \end{conj}

\begin{rem}
    We should emphasize that we could not discuss complex geometric properties of the both topological mirror candidates $(Y,\pi:Y \to \bbP^N)$ and $(\tilde{Y},\tilde{\pi}:\tilde{Y} \to \C^N)$ as we don't know how to glue complex structures. It means that the perverse Leray filtrations associated to $\pi$ and $\tilde{\pi}$ are not a priori well-defined in our case. Instead, we consider the potentially equivalent filtrations, called a general flag filtrations, whose description is purely topological. See Section \ref{sec:perverse} for the reviews. 
\end{rem}

  The motivational work for the appearance of the weight and perverse filtration in mirror symmetry, which we shall call \textit{mirror P=W phenomena}, goes back to the proposal of  Harder-Katzarkov-Przyjalkowski \cite{PWmirror} in the context of mirror symmetry of log Calabi-Yau's. For a given log-Calabi Yau manifold $U$ of complex dimension $n$, one can consider the mixed Hodge structure on the cohomology $H^*(U)$ that consists of Deligne's canonical weight filtration $W_\bullet$ and Hodge filtration $F^\bullet$. On the other hand, we have a canonical affinization map $\mathrm{Aff}:U \to \Spec H^0(U,\cO_U)$ and it provides the perverse Leray filtration $P_\bullet$. Since those filtrations are compatible with each other, we can define the perverse-mixed Hodge polynomial $PW_U$ by 
\begin{equation}
    PW_U(u,t,w,p):=\sum_{a,b,r,s}(\dim \Gr_F^a\Gr_{s+b}^W \Gr^P_{s+r}(H^s(U,\C)))u^at^sw^bp^r.
\end{equation}
\begin{conj}\label{conj:mirrorP=W}(Mirror P=W conjecture)\cite{PWmirror} 
     Assume that two $n$-dimensional log Calabi-Yau varieties $U$ and $U^\vee$ are mirror to each other. Then we have the following polynomial identity
     \begin{equation}\label{eq:mirrorP=W}
          PW_U(u^{-1}t^{-2}, t,p,w)u^nt^n=PW_{U^\vee}(u,t,w,p)
     \end{equation}
 \end{conj}   
 In case that $U$ has compactification $(X,D)$ where $X$ is a smooth (quasi-)Fano and $D$ is simple normal crossing anti-canonical divisor, the mirror P=W conjecture can be deduced from mirror symmetry for the pair $(X,D)$. Note that the choice of a pair $(X,D)$ corresponds to the choice of a hybrid LG potential $h:Y \to \C^N$ which plays a role of a proper affinization map. Then mirror symmetry expects that one could match the (part of) $E_1$-page of the spectral sequence for the weight filtration on $H^*(U)$ with the $E_1$-page of the spectral sequence for $G$-flag filtration (=perverse Leray filtration) associated to $h$ on $H^*(Y)$. Explicitly, this is an isomorphism of the $E_1$-pages
 \begin{equation}\label{eq: mirror PW intro}
     (\bigoplus_{p-q=a}\Gr_F^p{}^WE^{-l, p+q+l},d_1) \cong ({}^G E^{-l, n+a+l}, d_1^G)
 \end{equation}
 where both are known to degenerate at the $E_2$-page (See Section \ref{sec:weight},\ref{sec:perverse} for the notations). We say the mirror pair $(X, D)|(Y, \omega, h:Y \to \C^N)$ satisfies mirror P=W conjecture in a strong sense if the relation (\ref{eq: mirror PW intro}) holds. 
 
 \begin{thm}
    Suppose that each mirror pair  $(X_i, \cup_{j \neq i}X_{ij})|(Y_i, \omega_i, h_i:Y_i \to \C^N)$ satisfies mirror P=W conjecture in a strong sense. Then
\begin{enumerate}
    \item for $X$ and $Y$ as above, we have 
    \begin{equation}
    \bigoplus_{p-q=a}\Gr_F^p\Gr^{W_{\lim{}}}_{p+q}H^{p+q+l}(X) \cong \Gr^{P^\pi}_{n+a}H^{n+a+l}(Y)
    \end{equation}
    \item for $X_c$ and $\tilde{Y}$ as above, we have 
    \begin{equation}
        \bigoplus_{p-q=a}\Gr_F^p\Gr^{W}_{p+q}H^{p+q+l}(X_c) \cong \Gr^{P^{\tilde{\pi}}}_{n+a}H_c^{n+a+l}(\tilde{Y})
    \end{equation}
\end{enumerate}
\end{thm}

The main idea is to apply the gluing property (Proposition \ref{prop:Gluing Property in general}) of each hybrid LG potential $h_i:Y_i \to \C^N$ to describe the $E_1$-pages of the spectral sequences for  $P_\bullet^\pi$ and $P^{\tilde{\pi}}_\bullet$ in a way that they become isomorphic to those for $W_{\lim{}\bullet}$ and $W_\bullet$, respectively. One of the key lemmas is the Poincaré duality statement for hybrid LG models, which we will prove in Section \ref{sec:Appendix}.

\begin{thm}(Theorem \ref{thm:PD hybrid})(Poincaré duality)
    Let $(Y, h:Y \to \C^N)$ be a rank $N$ hybrid LG model. Then for $a \geq 0$, there is an isomorphism of cohomology groups 
    \begin{equation}
        H^a(Y, Y_{sm},\C) \cong H^{2n-a}(Y, Y_{sm},\C)^*
    \end{equation}
    where $n=\dim_\C Y$.
\end{thm}

\subsection{Semi-stable toric degeneration}

We provide an evidence for a mirror correspondence between a semi-stable degeneration of $X$ and a morphism from $Y$ to a projective variety for a Batyrev-Borisov mirror pair $(X,Y)$ \cite{BatyrevBorisovmirror}. We first consider a semi-stable degeneration of a smooth toric Fano variety $X_\Delta$ which is induced by a semi-stable partition $\Gamma$ of the polytope $\Delta$ (Definition \ref{def:semistable parition}). The degeneration fiber is the union of toric varieties associated to the maximal subpolytopes $\{\Delta_{(i)}|i=0, \dots, N\}$ in $\Delta$ for some $N$. This induces a type $(N+1)$ semi-stable degeneration of a general Calabi-Yau hypersuface $X$ of $X_\Delta$ whose degeneration fiber $X_c=\cup_{i=0}^N X_i$ is the simple normal crossing union of general hypersurfaces $X_i$ of $X_{\Delta_{(i)}}$, determined by $\Delta_{(i)}$. We will show that on the mirror side, the partition $\Gamma$ canonically induces a morphism $\pi:Y \to \bbP^N$ from a mirror dual Calabi-Yau $Y$. It follows from the construction that the deepest intersection of the components of the degeneration fiber is mirror to a generic fiber of $\pi$. Moreover, as the base of $\pi:Y \to \bbP^N$ comes with the toric chart, we obtain a natural candidate of a mirror hybrid LG model to $(X_i, \cup_{j \neq i}X_{ij})$. In other words, we take $\Delta_i:=\{|z_j| \leq |z_i||j \neq i\} \subset \bbP^N$ and set $Y_i:=\pi^{-1}(\Delta_i)$ and $h_i:=\pi|_{Y_i}$. 
    \begin{conj}\label{conj:hybrid LG model intro}(Conjecture \ref{Conj:hybrid LG model})
        For each $i$, the hybrid LG model $(Y_i,h_i: Y_i \to \Delta_i)$ is mirror to the pair $(X_i,\cup_{j \neq i}X_{ij})$. 
    \end{conj}

Conjecture \ref{conj:hybrid LG model intro} can be considered as the reverse construction of the topological gluing of hybrid LG models. We will discuss several issues for proving Conjecture \ref{conj:hybrid LG model intro} and leave it for future work. 

\subsection{Acknowledgement}
I would like to thank Nick Sheridan for
his consistent support and helpful conversation on the gluing arguments. I am also grateful to   Charles Doran, Andrew Harder and Alan Thompson for helpful discussion. The author gratefully acknowledges financial support from the Leverhulme Trust.
\section{Backgrounds}
In this section, we set up the notations and review basic concepts about mixed Hodge structures and perverse filtrations. 

\subsection{The weight filtration}\label{sec:weight}

We recall basic concepts about the mixed Hodge structures used in the main part of the article by following the exposition in \cite{peters2008mixed}. Let $U$ be a smooth quasi-projective variety over $\C$ and $(X,D)$ be a good compactification of $U$. Recall that a pair $(X,D)$ is called a good compactification of $U$ if $X$ is a smooth and compact variety and $D$ is a simple normal crossing divisor. Let $j:U \to X$ be a natural inclusion. Consider the logarithmic de Rham complex
\begin{equation*}
    \Omega^{\bullet}_X(\log D) \subset j_*\Omega^{\bullet}_U
\end{equation*}
Locally at $p \in D$ with an open neighborhood $V \subset X$ with coordinates $(z_1, \cdots, z_n)$ in which $D$ is given by $z_1\cdots z_k=0$, one can see
\begin{equation*}
    \begin{aligned}
    &\Omega^1_X(\log D)_p=\cO_{X,p}\frac{dz_1}{z_1} \oplus \cdots \oplus \cO_{X,p}\frac{dz_k}{z_k} \oplus \cO_{X,p}dz_{k+1} \oplus \cdots \oplus \cO_{X,p}dz_n\\
&\Omega^r_X(\log D)_p = \bigwedge^r\Omega^1_X(\log D)_p
    \end{aligned}
\end{equation*}

There are two natural filtrations on the logarithmic de Rham complex $(\Omega^{\bullet}_X(\log D),d)$:
    \begin{enumerate}
        \item (Hodge filtration) A decreasing filtration $F^\bullet$ on $\Omega^{\bullet}_X(\log D)$ defined by \begin{equation*}
            F^p\Omega^\bullet_X(\log D):= \Omega^{\geq p}_X(\log D) 
        \end{equation*}
        \item (Weight filtration) An increasing filtration $W_\bullet$ on $\Omega^{\bullet}_X(\log D)$ defined by 
        \[W_m\Omega^r_X(\log D):= \begin{cases}
        0 &  m<0\\
        \Omega^r_X(\log D) &  m \geq r\\
        \Omega^{r-m}_X\wedge \Omega^m_X(\log D) &  0 \leq m \leq r
        \end{cases}
        \]
    \end{enumerate}
\begin{thm}\cite[Theorem 4.2]{peters2008mixed}
    \begin{enumerate}
        \item The logarithmic de Rham complex $\Omega^{\bullet}_X(\log D)$ is quasi-isomorphic to $j_*\Omega^{\bullet}_U$:
        \begin{equation*}
         H^k(U;\C)=\mathbb{H}^k(X, \Omega^{\bullet}_X(\log D))
        \end{equation*}
        \item The decreasing filtration $F^\bullet$ on $\Omega^{\bullet}_X(\log D)$ induces the filtration in cohomology
        \begin{equation*}
          F^pH^k(U;\C)=\Image(\mathbb{H}^k(X, F^p\Omega^{\bullet}_X(\log D)) \to H^k(U;\C))
        \end{equation*}
        which is called the \textit{Hodge filtration} on $H^\bullet(U)$. Similarly, the increasing filtration $W_\bullet$ on $\Omega_X^\bullet( \log D)$ induces the filtration in cohomology
        \begin{equation*}
          W_mH^k(U;\C)=\Image(\mathbb{H}^k(X, W_{m-k}\Omega^{\bullet}_X(\log D)) \to H^k(U;\C))
        \end{equation*}
        which is called the \textit{weight filtration} on $H^\bullet(U)$. 
        \item The package $(\Omega^{\bullet}_X(\log D), W_{\bullet}, F^{\bullet})$ gives a $\C$-mixed Hodge structure on $H^k(U;\C)$.
    \end{enumerate}
\end{thm}

\begin{rem}
In general, the weight filtration can be defined over the field of rational numbers $\Q$ so that the $\Q$-mixed Hodge structures are considered. However, as we will mainly consider the filtrations on the cohomology group with complex coefficients, we do not emphasize the rational structures. 
\end{rem}

The key properties of these two filtrations are the degenerations of the associated spectral sequences. More precisely, we have 

\begin{proposition}\cite[Theorem 4.2, Proposition 4.3]{peters2008mixed}
    \begin{enumerate}
    \item The spectral sequence for $(\mathbb{H}(X, \Omega^{\bullet}_X(\log D)), F^{\bullet})$ whose $E_1$-page is given by 
    \begin{equation*}
         E_1^{p,q}=\mathbb{H}^{p+q}(X, \Gr_F^p\Omega^{\bullet}_X(\log D))
    \end{equation*}
    degenerates at $E_1$-page. Thus, we have
    \begin{equation*}
        \Gr_{F}^p\mathbb{H}^{p+q}(X, \Omega^{\bullet}_X(\log D))=\mathbb{H}^{p+q}(X, \Gr_F^p\Omega^{\bullet}_X(\log D)).
    \end{equation*}
    \item The spectral sequence for $(\mathbb{H}(X, \Omega^{\bullet}_X(\log D)), W_{\bullet})$ whose $E_1$-page is given by 
    \begin{equation*}
        E_1^{-m, k+m}=\mathbb{H}^k(X, \Gr^W_m\Omega^{\bullet}_X(\log D)) 
    \end{equation*}
    degenerates at $E_2$-page and the differential $d_1:E_1^{-m, k+m} \to E_1^{-m+1, k+m}$ is strictly compatible with the filtration $F_\bullet$. In other words, 
    \begin{equation*}
        E_2^{-m, k+m}=E_{\infty}^{-m, k+m}=\Gr^{W}_{m+k}\mathbb{H}^k(X, \Omega^{\bullet}_X(\log D)).
    \end{equation*}
    \end{enumerate}
\end{proposition}

For a given mixed Hodge structure $V=(V_\C, W_\bullet, F^\bullet)$ and $m \in \Z$, we define the $m$-th \textit{Tate twist} of $V$ by setting $V(m):=(V_\C(m), W(m)_\bullet, F(m)^\bullet)$ where $V_\C(m):=(2\pi i)^mV_\C$ and 
\begin{equation*}
    W(m)_k:=W_{k+2m} \quad F(m)^p:=F^{m+p}
\end{equation*}
for all $k$ and $p$.

In order to compute the mixed Hodge structure, we introduce the geometric description of the spectral sequence. Let $D$ be a simple normal crossing divisor with $N$ irreducible components $D_1, \dots, D_N$. For any index set $I \subset \{1, \dots, N\}$, we write $D_I=\cap_{i \in I}D_i$ for the intersection. We set $D(k)$ to be the disjoint union of $k$-tuple intersections of the components of $D$ and $D(0)$ to be $X$. Also, for $I=(i_1, \cdots, i_m)$ and $J=(i_1, \cdots, \Hat{i_j}, \cdots, i_m)$, there are inclusion maps 
\begin{align*}
    & \iota^I_J:D_I \hookrightarrow D_J\\
    \iota^m_j=\bigoplus_{|I|=m}& \iota_J^I:D(m) \hookrightarrow D(m-1)
\end{align*}
which induces canonical Gysin maps on the level of cohomology. Therefore, we have
\begin{equation}\label{eq:MV sign rule}
    \gamma_m=\oplus^m_{j=1}(-1)^{j-1}(\iota^m_j)_!:H^{k-m}(D(m))(-m) \rightarrow H^{k-m+2}(D(m-1))(-m+1) 
\end{equation}
where $(-)_!$ is the Gysin morphism. Here we call this sign convention by \textit{the Mayer-Vietoris sign rule} which is unique up to $\pm 1$. Under the residue map, this gives a geometric description of the differential $d_1:E_1^{-m, k+m} \to E_1^{-m+1, k+m}$ of the $E_1$-page of the spectral sequence for the weight filtration as follows:
\begin{proposition}\label{prop:CD d_1 spectral}\cite[Proposition 4.7]{peters2008mixed}
The following diagram is commutative
\begin{equation}\label{eq:E1 spectral weight}
    \begin{tikzcd}
        E_1^{-m, k+m} \arrow{r}{res_m} \arrow{d}{d_1} & H^{k-m}(D(m);\C)(-m) \arrow{d}{-r_m} \\
        E_1^{-m+1, k+m} \arrow{r}{res_{m-1}} & H^{k-m+2}(D(m-1);\C)(-m-1)
    \end{tikzcd}
\end{equation}
    where $res_m$ is the residue map for all $m \geq 0$.
\end{proposition}

Note that all morphisms in the diagram (\ref{eq:E1 spectral weight}) are compatible with Hodge filtration $F^\bullet$. This description provides several computational tools as well as functorial properties of mixed Hodge structures under geometric morphisms. For more details, we refer the reader to \cite{peters2008mixed}.

One can extend the above construction to the case when $U$ is singular. This can be done by using simplicial/cubical resolution of the singular variety $U$ and the notion of good compactification. We will not review this construction, but describe one particular case which we will deal with. We refer to \cite[Section 5]{peters2008mixed} for more details.

\begin{example}
    Let $D$ be a simple normal crossing variety with smooth components. The weight filtration $W_\bullet$ on $H^*(D,\C)$ has the following description. Consider the long exact sequences
    \begin{equation}\label{diag:spectral for D}
        0 \to \bigoplus_{|I|=1}H^j(D_I)\xrightarrow{d_0} \bigoplus_{|I|=2}H^j(D_I)\xrightarrow{d_1} \bigoplus_{|I|=3}H^j(D_I)\xrightarrow{d_2} \cdots 
    \end{equation}
    where $d_i$ is the alternating sum of the restriction map. Then we have 
    \begin{equation*}
        \Gr^W_{j}H^{i+j}(D)=\frac{\ker d_i}{\Image d_{i+1}}
    \end{equation*}
    In fact, the sequence (\ref{diag:spectral for D}) is the $E_1$-page of the spectral sequence for the weight filtration $W_\bullet$. It is also compatible with Hodge filtrations on each term to yield the Hodge filtration on $H^*(D)$. 
\end{example}

\subsection{The monodromy weight filtration}

Let $X$ be a smooth complex manifold and $\Delta$ be a unit disk. We consider a holomorphic map $f:X \to \Delta$ that is smooth over the punctured disk $\Delta^*:=\Delta \setminus \{0\}$. We also assume that $E:=f^{-1}(0)$ is a simple normal crossing divisor. Let $E_i$ be the components of $E$ and write
\begin{equation*}
    E_I=\bigcap_{i \in I} E_i, \qquad E(m)=\coprod _{|I|=m}E_I
\end{equation*}
as before. We present the de Rham theoretic description of the monodromy weight filtration on a generic fiber $X_t:=f^{-1}(t)$. 

Define the relative de Rham complex on $X$ with logarithmic poles along $E$:
\begin{equation*}
    \Omega^\bullet_{X/\Delta}(\log E):=\Omega_X^\bullet(\log E)/f^*(\Omega^1_\Delta (\log 0)) \wedge \Omega_X^{\bullet-1}(\log E)
\end{equation*}
By definition, this fits into the short exact sequence
\begin{equation*}
    0 \to f^*(\Omega^1_\Delta (\log 0)) \wedge \Omega_X^{\bullet-1}(\log E) \to \Omega_X^\bullet(\log E) \to \Omega^\bullet_{X/\Delta}(\log E) \to 0
\end{equation*}
By taking $\otimes \cO_E$ on the above sequence, we will have 
\begin{equation*}
    0\to \Omega^\bullet_{X/\Delta}(\log E)\otimes \cO_E [-1] \xrightarrow{\wedge dt/t} \Omega_X^\bullet(\log E) \otimes \cO_E \to \Omega^\bullet_{X/\Delta}(\log E)\otimes \cO_E \to 0
\end{equation*}
The connecting homomorphism induces the residue at 0 of the logarithmic extension of Gauss-Manin connection:
\begin{equation*}
    res_0(\nabla):\bbH^q(E, \Omega^\bullet_{X/\Delta}(\log E) \otimes \cO_E) \to \bbH^q(E, \Omega^\bullet_{X/\Delta}(\log E) \otimes \cO_E)
\end{equation*}
Note that the cohomology of the induced complex on $E$, $\Omega^\bullet_{X/\Delta}(\log E) \otimes \cO_E$ becomes isomorphic to the cohomology group $H^*(X_t)$. Also, the morphism $res_0(\nabla)$ is the same with recovers the monodromy action on $H^k(X_t)$.

Define the increasing filtration $W_\bullet$ on $\Omega_{X/\Delta}(\log E)\otimes \cO_E$ by 
\begin{equation*}
    W_k\Omega_{X/\Delta}(\log E)\otimes \cO_E:=\Image\{W_k\Omega_{X}(\log E) \to \Omega_{X/\Delta}(\log E)\otimes \cO_E\} 
\end{equation*}
and the decreasing filtration $F^\bullet$ by the simple truncation. To describe the monodromy weight filtration, we consider the resolution of $\Omega^\bullet_{X/\Delta}(\log E)\otimes \cO_E$ as follows. Define a tri-filtered double complex 
\begin{equation*}
    (A^{\bullet, \bullet}, d', d'', W_\bullet, W(M)_\bullet, F^\bullet)
\end{equation*}
on $E$ by 
\begin{equation*}
        \begin{aligned}
        &A^{p,q}=\frac{\Omega^{p+q+1}_X(\log E)}{W_p\Omega^{p+q+1}_Z(\log E)} \quad
        d'=(-) \wedge dt/t:A^{p,q} \to A^{p+1, q}\quad 
        d''=d_{dR}:A^{p,q} \to A^{p,q+1} \\
        & W_rA^{p,q}=\frac{W_{r+p+1}\Omega^{p+q+1}_X(\log E)}{W_p\Omega^{p+q+1}_Z(\log E)} \quad W(M)_rA^{p,q}=\frac{W_{r+2p+1}\Omega^{p+q+1}_X(\log E)}{W_p\Omega^{p+q+1}_Z(\log E)}\\
        & F^r A^{p,q}=\frac{F^r\Omega^{p+q+1}(\log E)}{W_p\Omega^{p+q+1}_Z(\log E)}
        \end{aligned}
    \end{equation*} 
    
We have the map
\begin{equation*}
    \begin{aligned}
        \mu: \Omega^q_{X/\Delta}(\log E)\otimes \cO_E & \to A^{0,q} \\  
            \omega & \mapsto (-1)^q(dt/t)\wedge \omega \mod{W_0}
    \end{aligned}
\end{equation*}
which defines a quasi-isomorphism of bi-filtered complexes 
\begin{equation*}
    \mu: (\Omega^\bullet_{X/\Delta}(\log E)\otimes \cO_E, W_\bullet, F^\bullet) \to (s(A^{\bullet, \bullet}), W_\bullet, F^\bullet)
\end{equation*}
where $s(A^{\bullet, \bullet})$ is the associated single complex. Consider the natural morphism $\nu:A^{p,q} \to A^{p+1, q-1}$ given by $\omega \mapsto \omega (\mod{W_{p+1}})$. As it commutes with both differentials $d'$ and $d''$, it induces the endomorphism of the associated simple complex $s(A^{\bullet,\bullet})$. Note that this sends $W(M)_r$ to $W(M)_{r-2}$ and $F^p$ to $F^{p-1}$.

\begin{thm} \cite[Theorem 11.21]{peters2008mixed}
    The following diagram is commutative:
    \begin{equation*}
    \begin{tikzcd}
        \bbH^q(E, \Omega^\bullet_{X/\Delta}(\log E)\otimes \cO_E) \arrow[r,"\mu"] \arrow[d, "res_0\nabla"] & \bbH^q(E, s(A^{\bullet, \bullet})) \arrow[d, "-\nu"] \\
        \bbH^q(E, \Omega^\bullet_{X/\Delta}(\log E)\otimes \cO_E) \arrow[r,"\mu"] & \bbH^q(E, s(A^{\bullet, \bullet}))
    \end{tikzcd}
    \end{equation*}
\end{thm}

By taking the residue map, we have 
\begin{equation*}
    \begin{aligned}
        \Gr^{W(M)}_r s(A^{\bullet, \bullet}) & \cong \bigoplus_{k \geq 0,-r}\Gr^W_{r+2k+1}\Omega_X^\bullet(\log E)[1] \\
        & \cong \bigoplus_{k \geq 0,-r}\Omega^\bullet_{E(r+2k+1)}[-r-2k]
    \end{aligned}
\end{equation*}
Therefore, the $E_1$ page of the spectral sequence for the monodromy weight filtration $W(M)_\bullet$ is given by 
\begin{equation*}
    E_1^{p,q}=\bigoplus_{k \geq 0,p}H^{q+2p-2k}(E(2k-p+1),\C)(p-k) \Longrightarrow H^*(X, \C)
\end{equation*}
which is compatible with Hodge filtrations. More explicitly, the $E_1$-page is given by the following diagram:

\begin{equation*}
    \begin{tikzcd}[column sep=0.6em]
         E_1^{-2, q} \arrow{r} & E_1^{-1, q}\arrow{r} & E_1^{0, q} \arrow{r}& E_1^{2, q}\\
         \bigoplus\limits_{k \geq 0} H^{q-2k-4}(E(2k+3)) & \bigoplus\limits_{k \geq 0} H^{q-2k-2}(E(2k+2)) & \bigoplus\limits_{k \geq 0}H^{q-2k}(E(2k+1))  &\bigoplus\limits_{k \geq 1}H^{q-2k+2}(E(2k)) \\
         H^{q-4}(E(3)) \arrow[r] \arrow[rd] \arrow[d, phantom, "\bigoplus"] & H^{q-2}(E(2)) \arrow[r] \arrow[rd] \arrow[d, phantom, "\bigoplus"]& H^q(E(1)) \arrow[rd]\arrow[d, phantom, "\bigoplus"] &   \\
         H^{q}(E(5)) \arrow[d, phantom, "\bigoplus"] \arrow[r] & H^{q-4}(E(4))\arrow[r] \arrow[rd] \arrow[d, phantom, "\bigoplus"] & H^{q-2}(E(3)) \arrow[r] \arrow[rd] \arrow[d, phantom, "\bigoplus"]& H^{q}(E(2))\\
         \vdots & \vdots & H^{q-4}(E(5))\arrow[r] & H^{q-2}(E(4)) \\
    \end{tikzcd}
\end{equation*}
where the horizontal arrows are (the alternating sum of) the Gysin morphisms while the anti-diagonal arrows are (the alternating sum of) the pullback homomorphism. If we write down two morpshisms by $G$ and $d$, respectively, the differential $d_1:E_1^{p,q} \to E_1^{p+1,q}$ is given by  $d_1=G+(-1)^pd$. 

\begin{thm}\cite[Theorem 11.22]{peters2008mixed}
    The spectral sequence for the filtration $W(M)_\bullet$ degenerates at $E_2$-page so that we have 
\begin{equation*}
    E_2^{p,q}=E_\infty^{p,q}=\Gr^{W(M)}_{q}H^{p+q}(X)
\end{equation*}
We will also denote the monodromy weight filtration $W(M)_\bullet$ by $W_{\lim{}\bullet}$.
\end{thm}

\subsection{The perverse filtration}\label{sec:perverse}
We briefly review the notion of perverse filtration \cite{Cataldotopologyofalgebraic} and its geometric description \cite{Cataldolef}.
\begin{definition}
    Let $Y$ be an algebraic variety (or scheme) with $D_c^b(Y)$ a derived category of constructible sheaves on $Y$. An object $K^{\bullet} \in D^b_c(Y)$ is called a \textit{perverse sheaf} if it satisfies following two dual conditions:
    \begin{enumerate}
        \item (Support Condition) $\dim \supp(\mathcal{H}^i(K^{\bullet})) \leq -i$
        \item (Cosupport Condition)
        $\dim \supp(\mathcal{H}^i(\mathbb{D}K^{\bullet})) \leq i$ where $\mathbb{D}:D^b_c(Y) \to D^b_c(Y)$ is a dualizing functor. 
    \end{enumerate}
\end{definition}

Verdier's dualizing functor on $D^b_c(Y)$ is defined as  $\mathbb{D}=\Hom_{\cO_Y}(-, p^{!}(\C_{pt}))$ where $p:Y \to pt$ is a trivial map. We call $p^{!}(\C_{pt})$ a dualizing complex of $Y$, and denote it by $\omega_Y$. In particular, if $Y$ is non-singular of complex dimension $n$, $\omega_Y=\C_Y[2n]$. Note that
the subcategory $\mathcal{P}(Y)$ of perverse sheaves on $Y$ is an abelian category. Also, the support and cosupport condition induces the so-called perverse $t$-structure $({}^\mathfrak{p}D^{b, \geq 0}_c(Y), {}^\mathfrak{p}D^{b, \leq 0}_c(Y))$ on $D^b_c(Y)$ whose heart is  $\mathcal{P}(Y)$. Explicitly, it is given by 
\begin{enumerate}
    \item $K^{\bullet} \in {}^\mathfrak{p}D^{b, \leq 0}_c(Y)$ if and only if $K$ satisfies the support condition. Also, ${}^\mathfrak{p}D^{b, \leq n}_c(Y):={}^\mathfrak{p} D^{b, \leq 0}_c(Y)[-n]$
    \item $K^{\bullet} \in {}^\mathfrak{p} D^{b, \geq 0}_c(Y)$ if and only if $K$ satisfies the cosupport condition. Also, ${}^\mathfrak{p}D^{b, \geq n}_c(Y):={}^\mathfrak{p} D^{b, \geq 0}_c(Y)[-n]$
\end{enumerate}
We denote ${}^\mathfrak{p}\tau_{\leq n}:D_c^b(Y) \to {}^\mathfrak{p} D^{b, \leq n}_c(Y)$ (resp. ${}^\mathfrak{p}\tau_{\geq n}:D_c^b(Y) \to {}^\mathfrak{p} D^{b, \geq n}_c(Y)$) the natural truncation functor. This induces \textit{perverse cohomology functors} ${}^\mathfrak{p}\mathcal{H}:D^b_c(Y) \to \mathcal{P}(Y)$ defined by ${}^\mathfrak{p}\mathcal{H}^k:={}^\mathfrak{p} \tau_{\leq 0} \circ {}^\mathfrak{p}\tau_{\geq 0} \circ [k]$. Applying the perverse truncation, one can define the perverse filtration on the hypercohomology of a constructible sheaf $\cK^\bullet$ on $Y$ as follows;
\begin{defin}
For $\cK^\bullet \in D^b_c(Y)$, the perverse filtration $P_\bullet$ on $\bbH^k(Y, \cK^\bullet)$ is defined to be
\begin{equation*}
    P_b\bbH^k(Y, \cK^\bullet):=\Image(\bbH^k(Y, {}^\mathfrak{p}\tau_{\leq b}\cK^\bullet) \to \bbH^k(Y, \cK^\bullet))
\end{equation*}
\end{defin}

Let $f:X \to Y$ be a morpshim of smooth varieties. Then we can define the perverse $f$-Leray filtration on the cohomology $H^\bullet(X, \C)$ by setting
\begin{equation*}
    P^f_lH^k(X,\C):=P^f_l\bbH^k(Y, Rf_*\C)
\end{equation*}
\begin{thm}\label{thm: perverse degenerate E_2}\cite[Corollary 14.41]{peters2008mixed}
    If $f$ is proper, then the perverse Leray spectral sequence degenerates at the $E_2$ page. In other words, we have 
    \begin{equation*}
        \Gr^{P^f}_l\bbH^{k}(X,\C)=E_2^{k-l,l}=\bbH^{k-l}(Y, {}^\mathfrak{p}\mathcal{H}^l(Rf_*\C))
    \end{equation*}
\end{thm}
We will provide geometric description of the perverse $f$-Leray filtration in case that the base space is either affine or quasi-projective. For this, we introduce some notations used in the next subsection. Let $f:X \to Y$ be a locally closed embedding. Then the restriction functor is given by $(-)|_X=Rf_!f^*$ on $D_c^b(Y)$, which is exact. If $f$ is closed, then we also have the right derived functor of sections with support in $X$, denoted by $R\Gamma_X(-)=Rf_*f^!$. 

\subsubsection{Geoemtric description}

In this subsection, we provide a geometric description of perverse filtrations. We follow the same convention for the indices of filtrations used in \cite{MirrorClemensSchmid}. Let's consider the following diagram of varieties
\begin{equation*}
    \begin{tikzcd}
         U \arrow[r, hook] \arrow[d, "\pi_U"] & Y \arrow[d, "\pi"] \\
         B_U \arrow[r, hook] & B 
    \end{tikzcd}
\end{equation*}
where 
\begin{itemize}
    \item $Y$ is a smooth complex projective variety of complex dimension $n$.
    \item $B$ is a complex projective variety of complex dimension $m$. We fix an embedding $B \hookrightarrow \bbP^N$. 
    \item $\pi:Y \to B$ is a projective and surjective morphism;
    \item $B_U$ is the affine subvariety of $B$ and $U:=\pi^{-1}(B_U)$. We write $\pi_U$ for the restriction of $\pi$ onto $U$.
\end{itemize}

Recall that there is a smooth projective variety $F(N,m)$ parametrizing $m$-flags $\mathfrak{F}=\{F_{-m} \subset \cdots \subset F_{-1}\}$ on $\bbP^N$ where $F_{-p}$ is a codimension $p$ linear subspace. A linear $m$-flag $\mathfrak{F}$ on $\bbP^N$ is general if it belongs to a suitable Zariski open subset of $F(N,m)$. Similarly, we say a pair of linear $m$-flags $(\mathfrak{F}_1, \mathfrak{F}_2)$ is general if it belongs to a suitable Zariski open subset of $F(N,m) \times F(N,m)$.

We fix a general pair of $m$-flags $({H}_\bullet, {L}_\bullet)$ on $\bbP^N$. Intersecting with $B$, it gives a pair of flags of subvarieties
$({B}_\bullet, {C}_\bullet)$ of $B$
\begin{equation*}
\begin{aligned}
    \emptyset=B_{-m-1} \subset B_{-m} \subset \dots \subset B_{-1} \subset B_0=B \\
    \emptyset=C_{-m-1} \subset C_{-m} \subset \dots \subset C_{-1} \subset C_0=B
\end{aligned}
\end{equation*}
where ${B}_\bullet:={H}_\bullet \cap B$ and ${C}_\bullet:={L}_\bullet \cap B$. We set $Y_\bullet=\pi^{-1}(B_\bullet)$ and $Z_\bullet=\pi^{-1}(C_\bullet)$. By following \cite{Cataldolef}\cite{MirrorClemensSchmid}, we define the following flag filtrations

\begin{defin}
\begin{enumerate}\label{def:flag filtrations}
    \item  The flag filtration $G^\bullet$ (of the first kind) on the cohomology of $U$ is the decreasing filtration 
    \begin{equation*}
        G^iH^k(U,\C):=\ker \{ \bbH^k(U, \C) \to \bbH^k(U, \C|_{U \cap Y_{i-1}}) \}
    \end{equation*}
    \item  The flag filtration $G^\bullet$ (of the second kind) on the compactly supported cohomology of $U$ is the decreasing filtration 
    \begin{equation*}
        G^jH_c^k(U,\C):=\Image\{ \bbH_{Z_{-j}\cap U,c}^k(U, \C) \to \bbH^k(U, \C) \}
    \end{equation*}
    \item  The $\delta$-flag filtration $\delta^\bullet$ on the cohomology of $Y$ is the decreasing filtration 
    \begin{equation*}
        \delta^p H^k(U,\C):=\Image\{ \bigoplus_{i+j=p} \bbH_{Z_{-j}}^k(Y, \C|_{Y-Y_{i-1}}) \to \bbH^k(Y, \C)\}
    \end{equation*}
\end{enumerate}
\end{defin}

Note that both two $G$-filtraions can be defined on other cohomology theories as well. We describe (the $E_1$-page of) the spectral sequence for each filtration.  
\begin{enumerate}
    \item The $E_1$-page of the spectral sequence for the flag filtration $G^\bullet$ (of the first kind) on $H^*(U)$ is given by 
    \begin{equation*}
        {}^GE_1^{p,q}=H^{p+q}(U \cap Y_p, U \cap Y_{p-1},\C) \Longrightarrow H^*(U,\C)
    \end{equation*}
    and the differential $d_1:{}^GE_1^{p,q} \to {}^GE_1^{p+1,q} $ is given by the connecting homomorphism of the long exact sequence of cohomology groups of the triple $(Y_{p+1}, Y_p, Y_{p-1})$. Furthermore, we have 
    \begin{equation*}
        {}^GE_\infty^{p,q} =\Gr_G^pH^{p+q}(U,\C)
    \end{equation*}
    \item The $E_1$-page of the spectral sequence for the flag filtration $G^\bullet$ (of the second kind) on $H_c^*(U,\C)$ is given by 
    \begin{equation*}
        {}^GE_{1}^{p,q}=H_{Z_{-p}\cap U-Z_{-p-1}\cap U,c}^{p+q}(U,\C) \Longrightarrow H_c^*(U,\C)
    \end{equation*}
    and the differential $d_1:{}^GE_{1}^{p,q} \to {}^GE_{1}^{p+1,q} $ is given by the connecting homomorphism of the long exact sequence of cohomology groups with supports $(Z_{-p}, Z_{-p-1}, Z_{-p-2})$. Furthermore, we have 
    \begin{equation*}
        {}^GE_{\infty}^{p,q} =\Gr_G^pH_c^{p+q}(U,\C)
    \end{equation*}
    \item The $E_1$-page of the spectral sequence for the $\delta$-flag filtration $\delta^\bullet$ on $H^*(Y)$ is given by 
    \begin{equation*}
    {}^\delta E_1^{p,q}=\bigoplus\limits_{i+j=p}H^{p+q}_{Z_{-j}-Z_{-j-1}}(Y, \C|_{Y_i-Y_{i-1}}) \Longrightarrow H^*(Y, \C)
\end{equation*}
    More explicitly, the $E_1$-page is given by the following diagram.
\begin{equation}\label{diag: delta E_1 page}
    \begin{tikzcd}[column sep=0.6em]
         E_1^{-2, q} \arrow{r} & E_1^{-1, q}\arrow{r} & E_1^{0, q} \arrow{r}& E_1^{1, q}\\
         H_{Z_0-Z_{-1}}^{q-2}( \C|_{Y_{-2}-Y_{-3}}) \arrow[d, phantom, "\bigoplus"] \arrow[r] \arrow[d, phantom, "\bigoplus"] \arrow[rd] & H_{Z_0-Z_{-1}}^{q-1}( \C|_{Y_{-1}-Y_{-2}}) \arrow[d, phantom, "\bigoplus"] \arrow[r] \arrow[rd]& H_{Z_0-Z_{-1}}^{q}( \C|_{Y_{0}-Y_{-1}}) \arrow[d, phantom, "\bigoplus"] \arrow[rd] &  &  \\
         H_{Z_{-1}-Z_{-2}}^{q-2}( \C|_{Y_{-3}-Y_{-4}}) \arrow[d, phantom, "\bigoplus"] \arrow[r]&  H_{Z_{-1}-Z_{-2}}^{q-1}( \C|_{Y_{-2}-Y_{-3}})\arrow[d, phantom, "\bigoplus"] \arrow[r] \arrow[rd]& H_{Z_{-1}-Z_{-2}}^{q}( \C|_{Y_{-1}-Y_{-2}})\arrow[d, phantom, "\bigoplus"]  \arrow[r] \arrow[rd]& H_{Z_{-1}-Z_{-2}}^{q+1}( \C|_{Y_{0}-Y_{-1}}) \arrow[d, phantom, "\bigoplus"] \\
         \vdots & \vdots & H_{Z_{-2}-Z_{-3}}^{q}( \C|_{Y_{-2}-Y_{-3}})\arrow[r] & H_{Z_{-2}-Z_{-3}}^{q+1}( \C|_{Y_{-1}-Y_{-2}}) \\
    \end{tikzcd}
\end{equation}
\begin{itemize}
    \item For fixed $j$, the anti-diagonal sequence is the same with the $E_1$-page of the spectral sequence for the $G$-filtration of the first kind on $H_{Z_{-j}}^*(Y)$ with respects to the induced flag $Z_{-j} \cap Y_\bullet$. Let's write $d_I$ for the differential.
    \item For fixed $i$, the horizontal sequence is the same with the $E_1$-page of the spectral sequence for the $G$-filtration of the second kind on $H^*(Y_i)$ with respects to the induced flag $Z_\bullet \cap Y_i$. Let's write $d_{II}$ for the differential.
    \item The differential $d_1:{}^\delta E_1^{p,q} \to {}^\delta E_1^{p+1,q}$ is given by $d_1=d_I+(-1)^pd_{II}$.
\end{itemize}

 Moreover, we have
\begin{equation*}
    {}^\delta E_\infty^{p,q}=\Gr_\delta^{p}H^{p+q}(Y,\C)
\end{equation*}
\end{enumerate}

\begin{thm}\label{thm:cataldo lef}\cite[Theorem 4.1.3 and 4.2.1]{Cataldolef}
    There are identification of the perverse and flag filtrations:
    \begin{enumerate}
        \item $P^{\pi_U}_lH^k(U)=G^{k-l}H^k(U)$ where $l$ starts from $k$ up to $k+m$.
        \item $P^{\pi_U}_lH_c^k(U)=G^{k-l}H_c^k(U)$ where $l$ starts from $k-m$ to $k$
        \item $P^\pi_lH^k(Y)=\delta^{k-l}H^k(U)$ where where $l$ starts from $k-m$ up to $k+m$.
    \end{enumerate}
\end{thm}

\begin{cor}
    The spectral sequences for all the flag filtrations in Definition \ref{def:flag filtrations} are degenerate at the $E_2$-page.
\end{cor}
\begin{proof}
Theorem \ref{thm:cataldo lef} implies that there are natural isomorphisms between two spectral sequences that induces the identity on the abutment. Since the morphisms $\pi$ and $\pi_U$ are proper, Theorem \ref{thm: perverse degenerate E_2} implies that the spectral sequence for the associated perverse filtration degenerates at the $E_1$-page. Note that the $E_2$-term of the Grothendieck spectral sequence used in Theorem \ref{thm: perverse degenerate E_2} is the same with $E_1$-term of the perverse filtration. Also, the spectral sequence for the shifted flag filtrations $G^{k-\bullet}$ and $\delta^{k-\bullet}$ is the shifted spectral sequence for $G^\bullet$ and $\delta^{\bullet}$, respectively. In other words, $E_1^{p,q}$ for the shifted filtration is the same with $E_2^{2p+q, -p}$ for the original filtration (see \cite{TheoriedeHodge2} \cite[Section 3.7]{Cataldolef}.) Therefore, we have the $E_2$-degeneration results for the flag filtrations $G^\bullet$ and $\delta^\bullet$.
\end{proof}

\section{Extended Fano/LG correspondence}\label{sec:extended Fano}
\subsection{Hybrid LG models}

We recall the notion of hybrid LG models introduced in \cite{mypapermirror}. Let's first introduce some notations. Let $h=(h_1, \dots, h_N):Y \to \C^N$ be a $N$-tuple of holomorphic functions and $(z_1,\dots, z_N)$ be a coordinate of the base $\C^N$. For each non-empty subset $I=\{i_1, \dots, i_l\} \subset \{1,\dots, N\}$, we write $h_I=(h_{i_1}, \dots, h_{i_l}):Y \to \C^{|I|}$ and the coordinate $(z_{i_1}, \dots, z_{i_l})$ for the base $\C^{|I|}$, which implicitly determines the natural inclusion $\C^{|I|} \subset \C^N$. 

\begin{defin}\label{def: weak hybrid LG model general}
	A \textit{hybrid Landau-Ginzburg (LG) model of rank $N$} is a triple 
	$(Y, \omega, h=(h_1, h_2, \dots, h_N):Y \to \C^N)$ where 
	\begin{enumerate}
		\item $(Y, \omega)$ is $n$-dimensional complex K\"ahler Calabi-Yau manifold with a K\"ahler form $\omega \in \Omega^2(Y)$;
		\item $h:Y \to \C^N$ is a proper surjective holomorphic map such that 
		\begin{enumerate}
		    \item (Local trivialization) there exists a constant $R>0$ such that for any non-empty subset $I \subset \{1, \dots, N\}$, the induced map $h_I:Y \to \C^{|I|}$ is a locally trivial symplectic fibration over the region $B_I:=\{|z_i| > R | i \in I\}$ with smooth fibers. Furthermore, over $B_I$ we have $v(h_j)=0$ for any horizontal vector field $v \in T^{h_I}Y$ associated to $h_I$ and $j \notin I$;
		    \item (Compatibility) for $I \subset J$,  such local trivializations are compatible under the natural inclusions $B_J \times \C^{N-|J|} \subset B_I \times \C^{N-|I|} \subset \C^N$. 
		\end{enumerate}
  	We call $h:Y \to \C^N$ a \textit{hybrid Landau-Ginzburg (LG) potential}. 
	\end{enumerate}
\end{defin}

Note that when $N=1$, it recovers the usual notion of LG models $(Y, \omega, h:Y \to \C)$ where $h$ becomes a locally trivial symplectic fibration with smooth fibers near the infinity. The second condition in Definition \ref{def: weak hybrid LG model general} generalizes the $N=1$ case, which controls the geometry of the local fibration $h$ near the infinity boundary on the base. For each non-empty subset $I \subset \{1, \dots, N\}$, let's write $Y_I$ for a generic fiber of $h_I:Y \to \C^{|I|}$ and $h_{Y_I}:Y_I \to \C^{N-|I|}$ for the restriction of $h$ into $Y_I$. Then the induced triple $(Y_I, \omega|_{Y_I}, h_{Y_I})$ can be regraded as a hybrid LG model of rank $N-|I|$. From this point of view, the condition $(a)$ in Definition \ref{def: weak hybrid LG model general} is rephrased as the condition that $h_I:Y \to \C^{|I|}$ is a locally trivialization of the induced hybrid LG models of rank $N-|I|$.

    Associated to the hybrid LG model $(Y, \omega, h:Y \to \C^N)$, we define the ordinary LG model to be a triple $(Y, \omega, \sw:=\Sigma \circ h:Y \to \C)$ where $\Sigma:\C^N \to \C$ is the summation map. To justify this terminology, we need to show that $\sw:Y \to \C$ is locally trivial near the infinity with smooth fibers. 

\begin{prop}\label{prop:Gluing Property in general}(Gluing Property)
    Let $(Y, \omega, h:Y \to \C^N)$ be a hybrid LG model and $H$ be a generic hyperplane in the base $\C^N$, which is not parallel to any coordinate hyperplanes. There exists an open cover $\{U_i\}_{i=1}^N$ of $H$ such that for any non-empty subset $I \subset \{1, \dots, N\}$, the induced map $h^{-1}(U_I) \to U_I$ is symplectic isotopic to the induced hybrid LG potential $h_{Y_I}:Y_I \to \C^{N-|I|}$ which is linear along the base. 
\end{prop}

\begin{proof}
    Take a hyperplane $H=\{a_1z_1 + \cdots + a_Nz_N=M\}$ where $a_i \neq 0$ for all $i$. By changing the coordinate $z_i \mapsto z_i/a_i$, we reduce to the case where $a_i=1$ for all $i$. We also further reduce to the case when $M$ is real due to the rotational symmetry. By generality, we take $M >NR$. First, note that $H \cap (\cap_{i=1}^N\{|z_i|\leq R\})=\emptyset$. Let $R_i=\{Re(z_i)>R\}$ and the simply-connected region
    \begin{equation*}
        U_i:=\{Re(z_i) > R\} \cap H=\{Re(z_1 + \cdots + \hat{z_i}+ \cdots + z_n)<M-R\} \cap H
    \end{equation*}
    for each $i$. Since $U_i \subset \{|z_i|>R\}$, one can project $U_i$ to the locus $\{z_i=2R\}$ inside the region $\{|z_i|<R\}$. The image of the projection is $V_i:=\{z_i=2R, Re(z_1+\cdots+\hat{z_i}+\cdots+z_N)<M-R\}$ which contains $\bigcap_{j \neq i}\{|z_j| \leq R\}$. Therefore, this projection identifies $h:h^{-1}(U_i) \to U_i$ with $h:h^{-1}(V_i) \to V_i$ due to the local triviality of the hybrid LG model. Moreover, the latter map is completed to $h_{Y_i}:Y_i \to \C^{N-1}$ by the inductive argument. In general, for each $I$, $U_I= \cap_{i \in I}U_i$ is non-empty and simply-connected. Since $U_I \subset \{|z_i|>R, i \in I\}$, one can apply the same argument to get the conclusion. 
\end{proof}

\begin{defin}
    Let $(Y, \omega, h:Y \to \C^N)$ be a hybrid LG model. We define the induced triple $(Y, \omega, \sw:Y \to \C)$ to be the ordinary LG model associated to the hybrid LG model $(Y, \omega, h)$ and denote a generic fiber of $\sw$ by $Y_{sm}$.
\end{defin}

On the cohomology level, Proposition \ref{prop:Gluing Property in general} implies that the cohomology group of $\pi^{-1}(H)$ are all (non-canonically) isomorphic to that of the normal crossing union of some $Y_i$'s. We will use this fact to study the perverse Leray filtration associated to $h:Y \to \C^N$ on $H^*(Y)$. 

Consider a general flag of hyperplanes in $\C^N$
\begin{equation*}
    \mathfrak{H}:0=H_{-N-1} \subset H_{-N} \subset \cdots \subset H_{-1} \subset H_0=\C^N
\end{equation*}
which is transversal to the discriminant locus of $h$ and each $H_{-l}$ is not parallel to any coordinate lines. We write $Y_{sm^{(l)}}$ for $h^{-1}(H_{-l})$ so that we have a general flag of subvarieties 
\begin{equation*}
    0 \subset Y_{sm^{(N)}} \subset \cdots \subset Y_{sm^{(1)}}\subset Y
\end{equation*}
which will compute the flag filtration $G^\bullet$ (equivalently, the  perverse Leray filtration $P^h_\bullet$ ) on $H^*(Y)$ (See Section \ref{sec:perverse}). In other words, the $E_1$-page of the spectral sequence is given by the long exact sequence 
\begin{equation}\label{eq:E_1-spectral perverse}
        H^{a-N}(Y_{sm^{(N)}}) \xrightarrow{d_1} H^{a-N+1}(Y_{sm^{(N-1)}}, Y_{sm^{(N)}}) \xrightarrow{d_1} \cdots \xrightarrow{d_1} H^{a-1}(Y_{sm^{(1)}}, Y_{sm^{(2)}}) \xrightarrow{d_1} H^a(Y, Y_{sm^{(1)}})
\end{equation}

We use the same notation in the proof of Proposition \ref{prop:Gluing Property in general}. Take open (simply-connected) regions $\{R_i \subset \C^N|i=1, \dots, N\}$ which induces the open covering of $H_{-1}$, $\{U_i:=R_i \cap H_{-1}|i=1, \dots, N\}$ that yields the gluing property. Let $V_i=\{z_i=const\}$ be the region that $U_i$ projects to. Due to the generality of the flag, we may assume that $H_{-2} \cap U_i \subset U_i$ projects to a hyperplane that is contained in $\cup_{j \neq i}R_j \cap V_i$ for all $i$. As $H_{-1}$ and $H_{-2}$ are not parallel to any coordinate lines, this can be done by scaling $M$ sufficiently large to place $H_{-2}$ far enough from the each coordinate line. Inductively, for each $k$, we could assume that the collection of regions $\{R_i \subset \C^N|i=1, \dots, N\}$ yields the gluing property for $H_{-k}$ in $V_I:=\cap_{i \in I}V_i$ for any $|I|=k-1$. Then the gluing property implies the following:
\begin{equation*}
    Y_{sm^{(k)}} \cap h^{-1}(U_I) \cong \begin{cases}
    Y_{I,sm^{(k-|I|)}} &|I| < k \\
    Y_I & |I| \geq k
    \end{cases}
\end{equation*}

\begin{lem}\label{lem:gluing description of flags}
For any $a \geq 0, k \geq 1$, the relative cohomology $H^a(Y_{sm^{(k)}}, Y_{sm^{(k+1)}})$ is isomorphic to $\bigoplus_{|I|=k}H^a(Y_I, Y_{I,sm})$.
\end{lem}

\begin{proof}
Take the (simply-connected) open region $\{R_i \subset \C^N|i=1, \dots, N\}$ and the induced cover $\{U_i=R_i \cap H_{-1}\}$ as above. When $k=1$, the Mayer-Vietoris argument with respects to the open cover $\{h^{-1}(U_i)\}$ and the gluing property implies that $H^a(Y_{sm^{(1)}}, Y_{sm^{(2)}}) \cong \bigoplus_{i=1}^N H^a(Y_i, Y_{i, sm^{(1)}})$ where $H^a(Y_i, Y_{i,sm^{(1)}})) \cong  H^a(Y_i,Y_{i,sm})$. In general, we apply the Mayer-Vietoris sequence to the cohomology group $H^a(Y_{sm^{(k)}}, Y_{sm^{(k+1)}})$ with the induced open cover by $R_i$'s. The $E_1$-page of the spectral sequence is given by 
\begin{equation*}
    \bigoplus_{|I|=1}H^a(Y_{I,sm^{(k-1)}},Y_{I,sm^{(k)}}) \xrightarrow{d_1} \dots \xrightarrow{d_1} \bigoplus_{|I|=k-1}H^a(Y_{I,sm^{(1)}},Y_{I,sm^{(2)}}) \xrightarrow{d_1}  \bigoplus_{|I|=k}H^a(Y_{I}, Y_{I,sm^{(1)}}) \rightarrow 0 ,
\end{equation*}
By induction, each direct summand is the direct sum of $H^a(Y_J, Y_{J,sm})$ for some $J$ with $|J|=k$. Then $d_1$ becomes the alternating sum of the identity morphisms where the signs are determined by the Mayer-Vietoris sign rule (\ref{eq:MV sign rule}). Then it follows from a simple combinatorial fact that this sequence is exact except at the first term, and $H^a(Y_{sm^{(k)}}, Y_{sm^{(k+1)}})=\ker(d_1) \cong \bigoplus_{|I|=k}H^a(Y_I, Y_{I,sm})$.
\end{proof}

For any $I \subset J$ with $|J|=|I|+1$, we write $\rho^{J}_I$ for the composition of morphisms
\begin{equation*}
    \rho^{J}_I:H^\bullet(Y_J, Y_{J, sm}) \hookrightarrow H^\bullet(Y_{I,sm}, Y_{I,sm^{(2)}}) \to  H^{\bullet+1}(Y_I, Y_{I, sm})
\end{equation*}
where the first one is given by Lemma \ref{lem:gluing description of flags} and the second one is the connecting homomorphism of the long exact sequence of cohomology groups of the triple 
$(Y_I, Y_{I,sm}, Y_{I, sm^{(2)}})$. Since we choose an open cover globally, Lemma \ref{lem:gluing description of flags} allows one to rewrite the $E_1$-page of the spectral sequence (\ref{eq:E_1-spectral perverse}) as follows:
\begin{equation*}
        H^{a-N}(Y_{sm^{(N)}}) \xrightarrow{d_1} \bigoplus_{|I|=N-1}H^{a-N+1}(Y_I, Y_{I,sm}) \xrightarrow{d_1} \cdots \xrightarrow{d_1} \bigoplus_{|I|=1}H^{a-1}(Y_I, Y_{I,sm}) \xrightarrow{d_1} H^a(Y, Y_{sm})
\end{equation*}
where the differential $d_1$ is the signed sum of the induced morphisms $\rho^{J}_I$'s that follows the Mayer-Vietoris sign rule (\ref{eq:MV sign rule})

For later use, we introduce the Poincaré dual of $\rho^J_I$. For a given hybrid LG model $(Y, h:Y \to \C^N)$ of rank $N$, we will show that there is a canonical isomorphism 
\begin{equation}\label{eq: PD hybrid}
    PD:H^a(Y, Y_{sm},\C) \xrightarrow{\cong} H^{2n-a}(Y, Y_{sm}, \C)^*
\end{equation}
for all $a \geq 0$ (Theorem \ref{thm:PD hybrid}). We define the morphism $(\rho^J_I)^\vee:H^\bullet(Y_I, Y_{I,sm}) \to H^{\bullet-1}(Y_J, Y_{J,sm})$ to be the composition $(\rho^J_I)^\vee=PD_J \circ (\rho^J_I)^* \circ PD_I^{-1}$ where $PD_I$ (resp. $PD_J$) is the same one in (\ref{eq: PD hybrid}) for the induced hybrid LG model $(Y_I, h_{Y_I})$ (resp. $(Y_J, h_{Y_J})$).

\subsection{Extended Fano/LG correspondence}

Let $X$ be a smooth (quasi-)Fano manifold and $D$ be an effective simple normal crossing anti-canonical divisor with $N$ components $D_1, D_2, \dots, D_N$. For any index set $I=\{i_1, i_2, \cdots, i_m\} \subset \{1, 2, \cdots, N\}$, we define
\begin{equation*}
    \begin{aligned}
    & D_I:=D_{i_1} \cap \cdots \cap D_{i_m} \\
    & D(I):=\Sigma_{j \notin I}D_I \cap D_j
    \end{aligned}
\end{equation*}
For example, if $I=\{1\}$, then $D_{\{1\}}=D_1$ and $D(\{1\})=(D_2 \cup \cdots \cup D_k) \cap D_1$. We also assume that all pairs $(D_I, D(I))$ are (quasi-)Fano. We also write the normal crossing union of $l$-th intersections by $D\{l\}=\sum_{|I|=l}D_I$ for all $l \geq 0$. 

\begin{defin} 
	A hybrid LG model $(Y, \omega, h:Y \to \C^N)$ is mirror to $(X, D)$ if it satisfies the following mirror relations:
	\begin{enumerate}
	    \item the associated ordinary Landau-Ginzburg model $(Y, \omega, \sw:Y \to \C)$ is mirror to $(X, D)$;
		\item for $i=1,2, \dots, N$, a hybrid Landau-Ginzburg model $(Y_i, \omega|_{Y_i}, h_{Y_i}:Y_i \to \C^{N-1})$ is mirror to $(D_i, D(\{i\}))$.
	\end{enumerate}
	Such a mirror pair is called a \textit{(quasi-)Fano mirror pair}, and we write it by $(X,D)|(Y, \omega, h)$.
\end{defin}
To elaborate the precise sense of the mirror relations, we introduce some notations. Let $\square$ be a cubical category whose objects are finite subsets of $\N$ and morphisms $Hom(I,J)$ consists of a single element if $I \subset J$ and else is empty. Given a category $\mathsf{C}$, we define a cubical object to be a contravariant functor $F:\square \to \mathsf{C}$, which we also call a cubical diagram of categories. For a cubical object $F$ and $I \subset \N$, we write
\begin{equation*}
    \begin{aligned}
        & F_I:=F(I)  \\
        & d_{IJ}:=F(I \to J): X_J \to X_I, \quad I \subset J
    \end{aligned}
\end{equation*}
We also define a morphism of cubical objects in an obvious way. We mainly consider the category of finite dimensional vector spaces over $\C$, denoted by $\mathsf{Vect_\C}$. 

First, on the $B$-side, consider the natural inclusions $\iota^J_I:D_J \hookrightarrow D_I$ for $I \subset J$. We consider the cubical object $\mathfrak{HH}_a(X,D)$ in $\mathsf{Vect}_\C$
\begin{equation*}
    \begin{aligned}
        &\mathfrak{HH}_a(X,D)_I=\bigoplus_{p-q=a}H^{p,q}(D_I) \\
        &\mathfrak{HH}_a(X,D)_{IJ}=\iota^J_{I!}:\bigoplus_{p-q=a}H^{p,q}(D_J) \to \bigoplus_{p-q=a}H^{p,q}(D_I)
    \end{aligned}
\end{equation*}
where $\iota_!$'s are the Gysin morphisms. We also consider the Poincaré dual of $\mathfrak{HH}_a(X,D)$, denoted by $\mathfrak{HH}^c_a(X,D)$ where
\begin{equation*}
    \begin{aligned}
        &\mathfrak{HH}^c_a(X,D)_I=\bigoplus_{p-q=a}H^{p,q}(D_I) \\
        &\mathfrak{HH}^c_a(X,D)_{IJ}=(\iota^J_{I})^*:\bigoplus_{p-q=a}H^{p,q}(D_I) \to \bigoplus_{p-q=a}H^{p,q}(D_J)
    \end{aligned}
\end{equation*}

On the A-side, let $(Y, \omega, h:Y \to \C^N)$ be a hybrid LG model of rank $N$ and $n=\dim_\C Y$. We consider the cubical object  $\mathfrak{HH}_a(Y,h) \in \mathsf{Vect}_\C$
\begin{equation*}
    \begin{aligned}
        &\mathfrak{HH}_a(Y,h)_I=H^{n+a-|I|}(Y_I, Y_{I,sm}) \\
        &\mathfrak{HH}_a(Y,h)_{IJ}=\rho^J_{I}:H^{n+a-|J|}(Y_J,Y_{J,sm}) \to H^{n+a-|I|}(Y_I,Y_{I,sm})
    \end{aligned}
\end{equation*}
for $n \leq a \leq n$. We also consider the Poincaré dual $\mathfrak{HH}^c_a(Y,h)$ where 
\begin{equation*}
    \begin{aligned}
        &\mathfrak{HH}^c_a(Y,h)_I=H^{n+a-|I|}(Y_I, Y_{I,sm}) \\
        &\mathfrak{HH}^c_a(Y,h)_{IJ}=(\rho^J_{I})^\vee:H^{n+a-|I|}(Y_I,Y_{I,sm}) \to H^{n+a-|J|}(Y_I,Y_{J,sm})
    \end{aligned}
\end{equation*}

\begin{conj}\label{conj: CMS general} Let $(X,D)|(Y,\omega, h:Y \to \C^N)$ be a (quasi-)Fano mirror pair. For $n \leq a \leq n$, there exists isomorphisms of cubical objects in $\mathsf{Vect}_\C$:
\begin{equation*}
        \mathfrak{HH}_a(X,D) \cong \mathfrak{HH}_a(Y,h) \qquad
        \mathfrak{HH}^c_a(X,D) \cong \mathfrak{HH}^c_a(Y,h)
\end{equation*}
\end{conj}

\begin{rem}
    Conjecture \ref{conj: CMS general} is motivated from the relative version of Homological mirror symmetry for (quasi-)Fano mirror pair \cite[Section 4.3]{mypapermirror}. In particular, this is expected to follow from applying Hochschild homology on the categorical statement. Also, it is expected that one of the above isomorphisms follow from the other via Poincaré duality. 
\end{rem}

\subsection{Line bundles/Monodromy correspondence}

Let $(X,D)$ be a quasi-Fano pair where $D$ is smooth and $(Y,\omega, \sw:Y \to \C)$ its mirror LG model. In this case, there is a mirror correspondence between the anti-canonical line bundle $-K_X$ and the monodromy $T$ of a generic fiber $\sw^{-1}(t)$ around the $\infty$. This correspondence is expected to be made precise on the categorical level via homological mirror symmetry conjecture. On the B-side, tensoring with $-K_X$, provides the auto-equivalence on the derived category of coherent sheaves on $X$, $D^b\Coh(X)$ as well as the one on $D^b\Coh(D)$ by the restriction. On the A-side, the monodromy operation $T$ induces auto-equivalences on the relevant Fukaya categories associated to $Y_{sm}$ and $\sw:Y \to \C$.

On the other hand, when $D$ has more than one component, one can ask more refined version of the above discussion. On the B-side, we have $N$ line bundles $\cO_X(D_i)$ for $i=1, \dots, N$, whose sum is the anti-canonical divisor $-K_X$. Each line bundle induces an auto-equivalence on $D^b\Coh(X)$ by taking a tensor product with itself. On the A-side, there are $N$ monodromy operations induced by taking a loop $T_i$ near the infinity on the base of $h:Y \to \C^N$
        \begin{equation}\label{eq:monoromies}
            T_i:=(t_1, \dots, t_{i-1}, e^{\sqrt{-1}\theta}t_i, t_{i+1}, \dots, t_N) \quad (0 \leq \theta \leq 2\pi)
        \end{equation}
        for a generic $(t_1, \dots, t_N) \in \C^N$ and $i=1, \dots, N$. We denote such operations by $\phi_{T_i}$. The monodromy $\phi_{T_i}$ induces not only the automorphism of a generic fiber $Y_i=h_i^{-1}(t)$ but also the automorphism of the induced fibration $h|_{Y_j}:Y_j \to \C^{N-1}$ for any $i,j$. This will play a key role in Section \ref{sec:mirror construction}. Moreover, note that the composition of $T_i$'s is the loop $T$ near the infinity on the base of $\sw:Y \to \C$. Each monodromy operation is expected to induce an auto-equivalence, denoted by $\phi_{T_i}$ as well, on the relevant Fukaya category of $(Y, \omega, h:Y \to \C^N)$.
    
\begin{Ansatz}\label{Ansatz:linebundle monodromy}
    There are correspondences between line bundle $\cO_X(D_i)$ and the monodromy $\phi_{T_i}$ for all $i=1, \dots, N$.
\end{Ansatz}

The main source of Ansatz \ref{Ansatz:linebundle monodromy} can be found in \cite{HanlonMonomdromy}\cite{HanlonHicksFunctorialityandhomological} where the mirror symmetry of smooth toric Fanos has been discussed. See also \cite[Section 4.3]{mypapermirror} for more details.

\section{Mirror construction for a smoothing of a semi-stable degeneration}\label{sec:mirror construction}
\subsection{Semi-stable degeneration}

Let $\fX$ be a complex connected analytic space  and $\Delta$ be a unit disc. A degeneration is a proper flat surjective map $\pi:\fX \to \Delta$ such that $\fX-\pi^{-1}(0)$ is smooth and the fiber $\fX_t$ is a compact Kähler manifold for every $t \neq 0$. The central fiber $\fX_0:\pi^{-1}(0)$ is called the degenerate fiber. Given the degeneration $\pi:\fX \to \Delta$ and for $t \neq 0$, we say that $\fX_t$ degenerates to $\fX_0$ or equivalently $\fX_0$ is smoothable to $\fX_t$. In particular, if the total space $\fX$ is smooth and the degenerate fiber $\fX_0$ is a simple normal crossing divisor of $\fX$, then the degeneration $\pi:\fX \to \Delta$ is called \textit{semi-stable}. We define \textit{type} of a semi-stable degeneration to be the dimension of the dual complex of the degenerate fiber. 

Due to Friedman, the semistability condition on the degeneration $\pi:\fX \to \Delta$ controls the behavior of the degenerate fiber in a way that the normal bundle of singular locus of $\fX_0$ in $\fX$ is trivial. This property is called \textit{ $d$-semistabilty}.

\begin{defin}\cite[Definition 1.13]{Friedmansmoothing}
    Let $X=\bigcup_{i=0}^N X_i$ be a normal crossing variety of pure dimension $n$ whose irreducible component is smooth. We define $X$ to be \textit{$d$-semistable} if 
    \begin{equation}\label{eq: d-semistability}
        \bigotimes_{i=0}^N I_{X_i}/I_{X_i}I_D \cong \cO_D
    \end{equation}
    where $D$ is the singular locus of $X$ and $I_{D}$ (resp. $I_{X_i}$) is the ideal sheaf of $I_{D}$ (resp. $I_{X_i})$.
\end{defin}

From now on, we specialize to the case where the degenerate fiber of the semi-stable degeneration of type $(N+1)$ consists of $N+1$ irreducible components. In this case, we have equivalent description of $d$-semistablity, which will be used in the mirror construction. Let's write $X_c=\bigcup_{i=0}^NX_i$ for the degenerate fiber $\fX_0$. For any $i,j \in \{0, \dots, N\}$, we write $X_{ij}$ for the intersection of $X_i$ and $X_j$ as a divisor of $X_i$. Since the degeneration $\pi:\fX \to \Delta$ is semi-stable, we have $X_c|_{X_i}=\fX_0|_{X_i} \cong \fX_t|_{X_i}=0$ for $t \neq 0$. It implies that in $\Pic(X_{ij})\cong \Pic(X_{ji})$, we have the following relation
\begin{equation*}
\begin{aligned}
    0&=\cO(X_0+\cdots+ X_N)|_{X_i}|_{X_{ij}}\\
    & =\cO(X_i)|_{X_{ij}}\otimes \cO(\sum_{j \neq i}\cO(X_j))|_{X_{ji}}
\end{aligned}
\end{equation*}
The right-hand side, denoted by $N(X_{ij})$, is called \textit{the normal class of $X_{ij}$}. A collection of all normal classes is the $\binom{N}{2}$-tuple
\begin{equation*}
    N_{X_c}:=(N(X_{ij})) \in \bigoplus_{i<j}\Pic(X_{ij})
\end{equation*}
and we call it a collection of normal classes of $X_c$. Note that the triviality of the collection of normal classes of $X_c$ implies (\ref{eq: d-semistability}). In our case, this is indeed equivalent. 

\begin{prop}
    Suppose that the normal crossing variety $X_c=\bigcup_{i=0}^N X_i$ introduced above is smoothable with a semi-stable degeneration of type $N+1$. Then the $d$-semistablity is equivalent to the triviality of the collection of normal classes of $X_c$. 
\end{prop}
\begin{proof}
    The same argument for the type III case \cite[Proposition 4.1]{Namhoonleetype3} applies to this case. 
\end{proof}

In general, the $d$-semistability condition is not sufficient to imply the smoothability with smooth total space. In case that $X_c$ is Calabi-Yau, which is of our main interest, this direction has been studied by Kawamata-Namikawa. 

\begin{thm}\cite[Theorem 4.2]{kawamatanamikawa}
    Let $X_c=\bigcup X_i$ be a compact Kähler normal crossing variety of dimension $n$ such that 
    \begin{enumerate}
        \item $X_c$ is $d$-semistable;
        \item its dualizing sheaf $\omega_{X_c}$ is trivial;
        \item $H^{n-2}(X_c, \cO_{X_c})=0$ and $H^{n-1}(X_i, \cO_{X_i})=0$ for all $i$.
    \end{enumerate}
    Then $X_c$ is smoothable to a Calabi-Yau $n$-fold $X$ with a smooth total space. 
\end{thm}
The following definition is motivated by the type III case \cite[Definition 2.1]{Namhoonleetype3}.
\begin{defin}
 Let $X$ be a Calabi-Yau projective normal crossing variety. $X$ is called \textit{d-semistable of type $(N+1)$} if there exists a type $(N+1)$ semi-stable degeneration $\phi:\fX \to \Delta$ whose central fiber $\fX_0$ is $X$.
\end{defin}

\begin{example}
Let $Q_5 \subset \bbP^4$ be a smooth quintic 3-fold. In the anti-canonincal linear system, it degenerates to a normal crossing union of smooth two hyperplanes $H_1$ and $H_2$ and a smooth cubic 3-fold $Q_3$. For simplicity, we denote it by $Z_c:=Z_1 \cup Z_2 \cup Z_3$ where $Z_1=H_1$, $Z_2=H_2$ and $Z_3=Q_3$. Note that the total space of such degeneration is singular so that one needs to modify $X_c$ to obtain a $d$-semistable degeneration. First, consider the intersection between a generic quintic 3-fold and $Z_c$. It becomes a union of three curves $C_1, C_2 \text{ and } C_3$ where $C_i$ lies in $Z_{jk}$ and $C_i \cap Z_{123}$ are all the same for $\{i,j,k\}=\{1,2,3\}$. For $\{i,j\}=\{1,2\}$, we take a blow up of $Z_i$ along $C_3$, denoted by $\pi_i:\Bl_{C_j}Z_i \to Z_i$. Let $E_i$ be an exceptional divisor and write $(-)'$ for the proper transformation of the subvariety $(-)$. While the proper transform $Z_{i3}'$ is isomorphic to $Z_{i3}$, $Z_{ij}'$ is the blow up of $Z_{ij}$ along $C_3 \cap Z_{ij}$. By construction $C_3'$ is disjoint from $Z_{123}'$. The last step is to blow up $\Bl_{C_3}Z_1$ along $C_3'$. If we write the resulting normal crossing variety as $X_c=X_1 \cup X_2 \cup X_3$, we have

\begin{equation}\label{eq: Quintic3-fold degeneration}
    \begin{aligned}
        & (X_1, X_{12} \cup X_{13}) \cong (\Bl_{C_3'}\Bl_{C_2}Z_1, Z_{12}' \cup Z_{13}) \\
        & (X_2, X_{21} \cup X_{23}) \cong (\Bl_{C_1}Z_2, Z_{12}' \cup Z_{13}) \\
        & (X_3, X_{31} \cup X_{32}) \cong (Z_3, Z_{12} \cup Z_{13}) \\
    \end{aligned}
\end{equation}
We will study how to construct a mirror of $Q_5$ from this degeneration data. 
\end{example}

\subsection{Mirror construction}

 Let $X_c=\bigcup_{i=0}^NX_i$ be a $d$-semistable Calabi-Yau $n$-fold of type $(N+1)$ and $X$ be a smoothing of $X_c$. For each $i$, the irreducible component $X_i$ is quasi-Fano with the canonically chosen anti-canonical divisor $\bigcup_{j \neq i}X_{ij}$. Suppose that the mirror hybrid LG model of the pair $(X_i, \cup_{j \neq i} X_{ij})$ is given by  $(Y_i, \omega_i, h_i=(h_{i0},\dots, h_{i\hat{i}}, \dots ,h_{iN}):Y_i \to \Delta^N)$. Here we shrink the base of the hybrid LG potential to the product of sufficiently large polydiscs $\Delta^N$. In this section, we will propose a topological construction of a mirror Calabi-Yau manifold of the smoothing $X$ by gluing the hybrid LG models $(Y_i, \omega_i, h_i)$. 
 
To perform the gluing, we require more topological conditions on the hybrid LG models. For $X_c$, two divisors $X_{ij} \subset X_i$ and $X_{ji} \subset X_j$ are topologically identified for $i \neq j$. This should be reflected on the mirror side by requiring that two induced hybrid LG models $(Y_{ij}:=h_{ij}^{-1}(t_j), h_i|_{Y_{ij}}:Y_{ij} \to \Delta^{N-1})$ and $(Y_{ji}:=h_{ji}^{-1}(t_i), h_j|_{Y_{ji}}:Y_{ji} \to \Delta^{N-1})$ are topologically the same for $t_i, t_j \in \partial \Delta$. Furthermore, one can enhance the topological identification by taking into account complex structures and symplectic structures. For instance, once the preferred choice of the topological identification $X_{ij}=X_{ji}$ has been made, the complex isomorphism between $X_{ij}$ and $X_{ji}$ is given by an element $f_{ij} \in \Aut(X_{ij})$ which is homotopic to the identity. Also these identifications are compatible to endow $X_c$ with a well-defined complex structure: For $i,j,k$, the composition of the restrictions $ f_{ki}|_{X_{kij}}\circ f_{jk}|_{X_{jki}} \circ f_{ij}|_{X_{ijk}} \in \Aut(X_{ijk})$ is homotopic to the identity. The mirror counterpart is the identification given by an element of symplectomorphisms $g_{ij} \in \Symp(Y_{ij}, \omega_i|_{Y_{ij}}, h_i|_{Y_{ij}})$ which preserves the hybrid LG potentials. Furthermore, for $i,j,k$, the composition of the restrictions $g_{ki}|_{(Y_{kij}, h_k)}\circ g_{jk}|_{(Y_{jki}, h_j)} \circ g_{ij}|_{(Y_{ijk}, h_i)} \in \Symp(Y_{ijk}, \omega_i|_{Y_{ijk}}, h_i|_{Y_{ijk}})$ is required to be homotopic to the identity.

\begin{Hypothesis}\label{Hypothesis}
    \begin{enumerate}
    \item For $i \neq j$, two induced hybrid LG models $(Y_{ij}:=h_{ij}^{-1}(t_j), h_i|_{Y_{ij}}:Y_{ij} \to \Delta^{N-1})$ and $(Y_{ji}:=h_{ji}^{-1}(t_i), h_j|_{Y_{ji}}:Y_{ji} \to \Delta^{N-1})$ are topologically the same for any $t_i, t_j \in \partial \Delta$. In particular, if symplectic structures are taken into account, this identification is given by a symplectomorphism $g_{ij} \in \Symp(Y_{ij}, \omega_i|_{Y_{ij}}, h_i|_{Y_{ij}})$, which is homotopic to the identity. 
    \item For $i \neq j \neq k$, the composition of the induced symplectomorphism $ g_{ki}|_{(Y_{kij}, h_k)}\circ g_{jk}|_{(Y_{jki}, h_j)} \circ g_{ij}|_{(Y_{ijk}, h_i)} \in \Symp(Y_{ijk}, \omega_i|_{Y_{ijk}}, h_i|_{Y_{ijk}})$ is homotopic to the identity.
\end{enumerate}
\end{Hypothesis}

Since there is already a global complex structure on $X_c$, without loss of generality, we may assume that all such gluing automorphisms are indeed the identity. This follows from perturbing complex and symplectic structures in the beginning.

Note that the identification on the bases $\Delta^N$ along the boundary components is modelled on the normal crossing union. For example, we can consider the base of the $i$-th hybrid LG model, denoted by $\Delta^N_{h_i}$, sits in $\C^{N+1}$ as 
\begin{equation*}
    \Delta^N_{h_i}\cong\{|z_j|\leq 1, z_i=t_i|j \neq i\}
\end{equation*}
for some $t_i$ with $|t_i|=1$. Thus we get a normal crossing union of $(Y_i, \omega_i)$ equipped with the induced map to a normal crossing union of the base $\Delta^N_{h_i}$ of each potential $h_i$. Moreover, topologically, we can further glue the bases $\Delta^N_{h_i}$ along the boundary components until the monodromies associated to $h_i$ come into this procedure. Then the resulting base becomes topologically equivalent to $\C^N$, hence we obtain a topological fibration $\tilde{\pi}:\tilde{Y} \to \C^N$. We will give more precise description of $\tilde{\pi}:\tilde{Y} \to \C^N$ later. 

From now on, we assume that the collection of the hybrid LG models $(Y_i, \omega_i, h_i:Y_i \to \Delta^N)$ satisfies Hypothesis \ref{Hypothesis}. Then we can interpret the vanishing of the normal classes of $X_c$ as the relation of the monodromies associated to hybrid LG models based on Ansatz \ref{Ansatz:linebundle monodromy}. Recall that $d$-semistability is equivalent to the triviality of normal class in $\Pic(X_{ij})$
\begin{equation}\label{eq:vanishing normal classeS}
\begin{aligned}
    0&=\cO(X_0+\cdots+ X_N)|_{X_i}|_{X_{ij}}\\
    & =\cO(X_i)|_{X_{ij}}\otimes \cO(\sum_{j \neq i}\cO(X_j))|_{X_{ji}}
\end{aligned}
\end{equation}
for any $i \neq j$. For each hybrid LG model $(Y_i, \omega_i, h_i:Y \to \Delta^N)$, we write monodromies induced by the loop along the $j$-th coordinate and the diagonal by $\phi_{T_{ij}} \text{ and } \phi_{T_i}$, respectively. Then the relation (\ref{eq:vanishing normal classeS}) corresponds to 
\begin{equation}\label{eq:gluing monodromy}
    \phi_{T_{ij}} \circ \phi_{T_{j}} = \mathrm{Id} \in \Symp(Y_{ij}, \omega_{i}|_{Y_{ij}}, h_{i}|_{Y_{ij}})
\end{equation}

In other words, we have the following correspondence of monodromies
\begin{equation*}
    \begin{tikzcd}
        Y_i \arrow[d, "h_i"] & Y_j \arrow[d, "h_j"] \\
        \Delta^{N} & \Delta^{N} \\
        (t_{i0}, \dots,e^{\sqrt{-1}\theta}t_{ij},\dots t_{iN}) \arrow[r, leftrightarrow] & e^{-\sqrt{-1}\theta}(t_{j0}, \dots, t_{j\hat{j}}, \dots, t_{jN}) 
    \end{tikzcd}
\end{equation*}

\begin{prop}\label{prop:Gluing LG models}
    Suppose that the hybrid LG models $\{(Y_i, \omega_i, h_i:Y_i \to \Delta^N)|i=0, \dots, N\}$ introduced above satisfy Hypothesis \ref{Hypothesis} and the relation (\ref{eq:gluing monodromy}). Then they can be glued to yield a symplectic fibration $\pi:Y \to \bbP^{N}$. 
\end{prop}

\begin{proof}
Let $[u_0:\cdots:u_N]$ be homogeneous coordinates on $\bbP^N$. Consider the closed subsets $\Delta_i \subset \bbP^N$
 \begin{equation*}
    \Delta_i=\{ |z_j| \leq |z_i| \text{ for } j \neq i\}
\end{equation*}
for $i=0, \dots, N$. First, note that $\cup_{i=0}^N \Delta_i=\bbP^N$. Since $\Delta_i \subset U_i=\{u_i \neq 0\}$, we have $\Delta_i \cong \Delta^N$ and any $k$-th intersection of $\Delta_i$'s is $(S^1)^k \times \Delta^{N-k}$. Due to Hypothesis \ref{Hypothesis}, we may identify the base of $h_i:Y \to \Delta^N$ with $\Delta_i$ for all $i$'s. Then it is straightforward to see that the relation (\ref{eq:gluing monodromy}) is the same with the induced chart maps for $\Delta_i$'s. 
\end{proof} 

 We keep the notation used in the proof of Proposition \ref{prop:Gluing LG models}. Consider the moment map $\mu:\bbP^{N}\to \bbR^{N}$ which is given by 
\begin{equation*}
    \begin{aligned}
        \mu:\bbP^{N}&\to \bbR^{N} \\
        [u_0:\cdots:u_N] &\mapsto \left(\frac{|u_0|}{\sum_{i=0}^N|u_i|}, \cdots, \frac{|u_{N-1}|}{\sum_{i=0}^N|u_i|}\right)
    \end{aligned}
\end{equation*}
Note that the image $\Image(\mu) \in \R^N$ is the standard $N$-simplex $\Delta$. Also, a fiber of $k$-dimensional face $\sigma$ is $(S^{1})^{k}$. We consider the dual spine $\Pi^N$ in $\Delta$, which is the subcomplex of the first barycentric subdivision of $\Delta$ spanned by the 0-skeleton of the first barycentric subdivision minus the 0-skeleton of $\Delta$. Decomposing along $\Pi^N$ gives $\Delta$ a natural cubical structure where the top-dimensional cubes pull back to $N$-th product of the disc $D^N$ in $\bbP^N$. These are exactly the polydiscs $\Delta_0, \dots, \Delta_N$ we have introduced. Each $k$-dimensional cube pull backs to $(S^1)^{N-k} \times D^k$ that is the $k$-th intersection of $\Delta_i$'s. 

Now we can see how $\tilde{\pi}:\tilde{Y} \to \C^N$ sits in the fibration $\pi:Y \to \bbP^N$ more rigorously. Take a standard open cover $\{V'_i|i=0, \dots, N\}$ of the $N$-simplex $\Delta$ where each $k$-th intersection of $V'_i$'s contains only one of codimension $k$ cubes in $\Delta$ for all $k$. We write $V_i$ for the preimage of $V'_i$ under $\mu$. Suppose we remove the image of a generic $\bbP^{N-1} \subset \bbP^N$ near the spine $\Pi^N$. The overlap $V_{ij}$ is now diffeomorphic to $\Delta^{N-1} \times (S^1 \times [0,1] -\{pt\})$ hence not contracts to $\Delta^{N-1} \times S^1$. However, if one instead removes a small closed neighborhood $N_\epsilon(\bbP^{N-1})$ of $\bbP^{N-1}$ in $\bbP^N$ and shrink $V_i$'s if necessary, then $V_{ij}$ becomes diffeomorphic to $\Delta^{N-1}$ (See Figure \ref{fig:figure1} for $N=3$). Since $V_i$ contracts to the base of the hybrid LG model $(Y_i, \omega_i, h_i)$, the induced symplectic fibration $\pi:{\pi^{-1}(\bbP^N \setminus N_\epsilon(\bbP^{N-1}))} \to \bbP^N \setminus N_\epsilon(\bbP^{N-1})) \cong \C^N$ can be seen as $\tilde{\pi}:\tilde{Y} \to \C^N$. In other words, Proposition \ref{prop:Gluing LG models} is equivalent to saying that $\tilde{\pi}:\tilde{Y} \to \C^N$ is compactifiable (to $\pi:Y \to \bbP^N$) if the condition (\ref{eq:gluing monodromy}) holds. 

\begin{center}
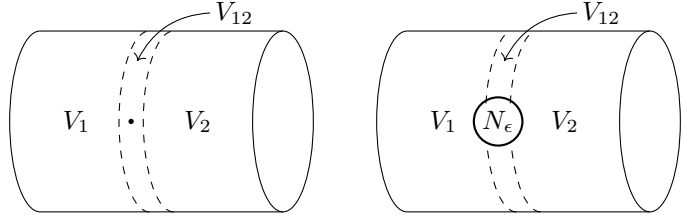

    \begin{tikzpicture}[scale = 0.8]
 \draw[dashed] (2.2,0) arc (-90:90:-0.5 and 1.5);% line 1
  \draw[dashed] (1.8,0) arc (-90:90:-0.5 and 1.5);% line 2

  \draw (0,0) -- (4,0);% bottom line
  \draw (0,3) -- (4,3);% top line
  \draw (0,0) arc (270:90:0.5 and 1.5);% left half of the left ellipse
  \draw (4,1.5) ellipse (0.5 and 1.5);% right ellipse
  \draw (0.6,1.5) node {$V_1$};
  \draw (2.6,1.5) node {$V_2$};
   \node at (1.5,1.5) [circle,inner sep=0.7pt,fill=black] {};
   \draw[->, bend right] (2.8, 3.3) to (1.6, 2.5);
   \node at (3.2, 3.3) {$V_{12}$};
  \end{tikzpicture}
  \qquad
\begin{tikzpicture}[scale = 0.8]
  \draw[dashed] (2.2,0) arc (-90:-15:-0.5 and 1.5);% line 1
  \draw[dashed] (1.3,
  1.8) arc (10:75:-0.5 and 1.5);
  \draw[dashed] (1.8,0) arc (-90:-15:-0.5 and 1.5);% line 2
  \draw[dashed] (1.7,1.8) arc (10:75:-0.5 and 1.5);%
  \draw (0,0) -- (4,0);% bottom line
  \draw (0,3) -- (4,3);% top line
  \draw (0,0) arc (270:90:0.5 and 1.5);% left half of the left ellipse
  \draw (4,1.5) ellipse (0.5 and 1.5);% right ellipse
  \draw (0.6,1.5) node {$V_1$};
  \draw (2.6,1.5) node {$V_2$};
  \draw[thick] (1.5,1.5) circle (4mm) node {$N_\epsilon$};
  \draw[->, bend right] (2.8, 3.3) to (1.6, 2.5);
   \node at (3.2, 3.3) {$V_{12}$};
  \end{tikzpicture}
  \captionof{figure}{Description of the intersection of $V_1$ and $V_2$ in $\bbP^2$}\label{fig:figure1}
\end{center}

\begin{rem}
 In general, there is a significant difference between gluing polydisks and standard open charts of $\bbP^N$. This is because the former procedure encodes information about singularities of each hybrid LG model while the latter is too rigid to do so. 
\end{rem}

\begin{thm}\label{thm: topological mirror symmetry}
    Let $X_c=\bigcup_{i=0}^NX_i$ be a $d$-semistable Calabi-Yau $n$-fold of type $(N+1)$ and $X$ be its smoothing. Suppose that we have hybrid LG models $(Y_i, \omega_i, h_i:Y_i \to \Delta^N)$ mirror to $(X_i, \cup_{j \neq i}X_{ij}$ that satisfies Hypothesis \ref{Hypothesis} and the relation (\ref{eq:gluing monodromy}). Let $\tilde{Y}$ and $Y$ be glued symplectic manifolds constructed above. Then
    \begin{enumerate}
        \item $Y$ is topological mirror to $X$. In other words, $e(Y)=(-1)^ne(X)$
        \item $\tilde{Y}$ is topological mirror to $X_c$. In other words, $e(\tilde{Y})=(-1)^ne(X)$
    \end{enumerate}
    where $e(-)$ is the Euler characteristic.
\end{thm}

\begin{lem}\label{lem:MV-Euler Characteristic}
   Let $h:Y \to \C^N$ be a hybrid LG model. Then $e(Y_{sm})=\sum_{|I|=1}^N(-1)^{|I|-1}e(Y_I)$. 
\end{lem}

\begin{proof}
This follows from the gluing property (Proposition \ref{prop:Gluing Property in general}). It implies that there exists an open cover $\{U_i|i=1, \dots, N\}$ of $Y_{sm}$ such that the induced fibration $h|_{U_I}:U_I \to \C^{N-|I|}$ is isotopic to $h|_{Y_I}: Y_I \to \C^{N-|I|}$. The conclusion follows from the Mayer-Vietoris argument. 
\end{proof}

\begin{proof}[Proof of Theorem \ref{thm: topological mirror symmetry}]
    Both items $(1)$ and $(2)$ are proven by the Mayer-Vietoris argument. We prove the item (2) first. 
    \begin{etaremune}
        \item By the Mayer-Vietoris sequence, we have
        \begin{equation*}
        \begin{aligned}
        e(X_c)&=\sum_{|I|=1}^{N+1}(-1)^{|I|-1}e(X_I) \\
         &=\sum_{|I|=1}^{N+1}(-1)^{|I|-1}(-1)^{n-|I|+1}e(Y_I, Y_{I,sm}) \\
         &=(-1)^{n}\sum_{|I|=1}^{N+1}e(Y_I, Y_{I,sm})
        \end{aligned}    
        \end{equation*}
        Now, it is enough to show that $e(\tilde{Y})=\sum_{|I|=1}^{N+1}e(Y_I, Y_{I,sm})$. Take an open cover $\{U_i|i=0, \dots, N\}$ of the base of $\tilde{\pi}:\tilde{Y} \to \C^N$ such that $\tilde{\pi}^{-1}(V_I)$ contracts to $Y_{I}$ for all $I$. By applying the Mayer-Vietoris argument, we have 
        \begin{equation*}
          e(\tilde{Y})=\sum_{|I|=1}^{N+1} (-1)^{|I|-1}e(Y_I)
        \end{equation*}
        
        Since $e(Y_I)=e(Y_I, Y_{I,sm}) + e(Y_{I,sm})$, Lemma \ref{lem:MV-Euler Characteristic} implies that 
        \begin{equation*}
        \begin{aligned}
        e(\tilde{Y})&=\sum_{|I|=1}^{N+1}(-1)^{|I|-1}(e(Y_I, Y_{I,sm})+e(Y_{I,sm}))\\
        &=\sum_{|I|=1}^{N+1}(-1)^{|I|-1}\left(e(Y_I, Y_{I,sm})+\sum_{|J|=1, J \cap I  = \emptyset}^{N+1-|I|}(-1)^{|J|-1}e(Y_{I,J})\right)
        \end{aligned}    
        \end{equation*}
        Here $Y_{I,J}$ is the same with $Y_{I \cup J}$ by Hypothesis \ref{Hypothesis} (1), but we use different notation to emphasize the Mayer-Vietoris procedure. By rewriting $e(Y_{I,J})=e(Y_{I,J},Y_{I,J,sm})+e(Y_{I,J,sm})$, we can iteratively apply Lemma \ref{lem:MV-Euler Characteristic}. Then we get $e(\tilde{Y})=\sum_{|I|=1}^{N+1}e(Y_I, Y_{I,sm})$ by taking the resummation.
        
        \item It follows from the same method in \cite[Proposition 3.2]{Namhoonleetype3} that the Euler characteristic of the smoothing manifold $X$ is given by 
        \begin{equation*}
          e(X)=\sum_{|I|=1}^{N+1}(-1)^{|I|}|I|e(X_I)
        \end{equation*}
        Similar to $(1)$, we have $e(X)=(-1)^{n}\sum_{|I|=1}^{N+1}|I|e(Y_I, Y_{I,sm})$. Take an open cover $\{U_i|i=0, \dots, N\}$ of the base of $h:Y \to \bbP^N$ as before. Note that for any $I$, the intersection $U_I$ contracts to $(S^1)^{|I|-1} \times \Delta^{N+1-|I|}$. Applying the Mayer-Vietoris argument with respects to the induced open cover $\{h^{-1}(U_i)\}$, one can see that the Euler characteristic $e(h^{-1}(U_I))$ vanishes for $|I| >1$ because $h^{-1}(U_I)$ contracts to a fiber bundle over $(S^1)^{|I|-1}$ with a fiber $Y_I$. Therefore, 
        \begin{equation*}
          e(Y)=\sum_{|I|=1}e(Y_I)=\sum_{|I|=1}e(Y_I, Y_{I,sm})+e(Y_{I,sm})
        \end{equation*}
        By applying Lemma \ref{lem:MV-Euler Characteristic} interatively, we get the conclusion.
    \end{etaremune}
\end{proof}

\section{Mirror P=W conjecture}
We keep the notation used in the previous section. Let $X_c=\bigcup_{i=0}^NX_i$ be a $d$-semistable Kähler Calabi-Yau $n$-fold of type $(N+1)$ and $X$ be its smoothing. We've introduced the topological construction of their mirror objects with additional symplectic fibration structure, $(\tilde{Y}, \tilde{\pi}:\tilde{Y} \to \C^N)$ and $(Y, \pi:Y \to \bbP^N)$, respectively. In this section, we discuss a refined version of Theorem \ref{thm: topological mirror symmetry}. In the degeneration picture, we have two filtrations on the cohomology groups: Deligne's canonical weight filtration $W_\bullet$ on $H^*(X_c)$ and the monodromy weight filtration $W_{\lim{}\bullet}$ on $H^*(X)$. The corresponding mirror filtrations are conjectured to be the perverse Leray filtration (equivalently, the $G$-flag filtration) on $H^*(\tilde{Y})$ associated to $\tilde{h}$ and the perverse filtration (equivalently, the $\delta$-flag filtration) on $H^*(Y)$ associated to $h$, respectively. 

\begin{center}
\begin{tabular}{ c | c }
 Degeneration (B-side) & Fibration (A-side)\\
 \hline 
 Deligne's canonical weight filtration $W_\bullet$ & Perverse filtration $P^{\tilde{\pi}}_\bullet$ associated to $\tilde{\pi}$ \\  
 Monodromy weight filtration $W_{\lim{}\bullet}$ & Perverse filtration $P^\pi_\bullet$ associated to $\pi$
\end{tabular}
\end{center}

\begin{rem}
Since the gluing construction (Proposition \ref{prop:Gluing LG models}) does not provide complex structures on the mirror objects, we could not discuss complex geometric properties of the mirrors, especially the perverse Leray filtrations. However, in our case, the main theorem of \cite{Cataldolef} shows that such filtrations are the same with general flag filtrations which has purely topological origin, as reviewed in \ref{sec:perverse}. The perverse Leray filtrations used in this section are indeed corresponding flag filtrations. 
\end{rem}

 \begin{thm}\label{thm:P=W main}
    Suppose that each mirror pair  $(X_i, \cup_{j \neq i}X_{ij})|(Y_i, \omega_i, h_i:Y_i \to \C^N)$ satisfies Conjecture \ref{conj: CMS general}.
\begin{enumerate}
    \item For $X$ and $Y$, we have 
    \begin{equation*}
    \bigoplus_{p-q=a}\Gr_F^p\Gr^{W_{lim}}_{p+q}H^{p+q+l}(X) \cong \Gr^{P^\pi}_{n+a}H^{n+a+l}(Y)
    \end{equation*}
    where $-N \leq l \leq N$.
    \item For $X_c$ and $\tilde{Y}$, we have
    \begin{equation*}
        \bigoplus_{p-q=a}\Gr_F^p\Gr^{W}_{p+q}H^{p+q+l}(X_c) \cong \Gr^{P^{\tilde{\pi}}}_{n+a}H_c^{n+a+l}(\tilde{Y})
    \end{equation*}
    where $0 \leq l \leq N$.
\end{enumerate}
\end{thm}

\begin{proof}[Proof of Theorem\ref{thm:P=W main}-(1)]

Let $\pi:Y \to \bbP^N$ be a gluing of ${N+1}$ hybrid LG models $(Y_i,h_i:Y_i \to \C^N)$. Take a standard open cover $\{V_s\}$ of $\bbP^N$ as explained in the discussion after Proposition \ref{prop:Gluing LG models} such that the induced fibration $\pi^{-1}(V_s) \to V_s$ contracts to $h_s:Y_s \to \Delta^N$. We consider a general linear flag in $\bbP^N$ 
\begin{equation*}
    \mathfrak{H}:0=H_{-N-1} \subset H_{-N} \subset \cdots \subset H_{-1} \subset H_0=\bbP^N
\end{equation*}
which satisfies several properties:
\begin{enumerate}
    \item $\mathfrak{H}$ intersects transversally with the discriminant locus of $\pi$;
    \item the induced flag $\mathfrak{H} \cap V_s$ is not parallel to any coordinate lines of the base of $h_s:Y_s \to \Delta^N$. 
\end{enumerate}

We also consider a pair of such general linear flags  $(H_\bullet, L_\bullet)$ of $\bbP^N$
\begin{equation*}
\begin{aligned}
    \mathfrak{H}:H_{-n} \subset \dots \subset H_{-1} \subset H_0=\bbP^N \\
    \mathfrak{L}:L_{-n} \subset \dots \subset L_{-1} \subset H_0=\bbP^N
\end{aligned}
\end{equation*}
Due to the generality, there exists a collection of sufficiently small $\epsilon_j>0$ which yields isomorphisms of pairs
\begin{equation}\label{eq:tubular pair}
    ((H_{-j}-H_{-j-1})\cap L_i, (H_{-j}-H_{-j-1}) \cap L_{i-1}) \cong ((H_{-j}-N_{\epsilon_j}(H_{-j-1}))\cap L_i, (H_{-j}-N_{\epsilon_j}(H_{-j-1})) \cap L_{i-1}) 
\end{equation}
for $i=-N, \dots, 0$ and $j=0, \dots, N$. Given this pair of flags, we will modify the cover $\{V_s\}$ in a way that the induced cover of $\bbP^N-N_{\epsilon_0}(H_{-1})$, denoted by $\{V_s^{(0)}\}$, satisfies the following properties:
\begin{enumerate}
    \item for non-empty index set $I$, $V^{(0)}_I \cong \Delta^{N+1-|I|}$;
    \item the induced the open regions $\{V^{(0)}_{st}|t \neq s\}$ yields the gluing property of the induced flag $V_s^{(0)} \cap \mathfrak{L}$ in the sense of the discussion before Lemma \ref{lem:gluing description of flags}.
\end{enumerate}
The first condition is one explained in the discussion after Proposition \ref{prop:Gluing LG models}. The second condition can be obtained by rescaling each hybrid LG potential $h_s:Y_s \to \Delta^N$ before the gluing. Next, we further modify the open cover $\{V_s\}$ in a way that the induced cover of $H_{-1} \setminus N_{\epsilon_1}(H_{-2})$, denoted by $\{V^{(1)}_s\}$ satisfies the following properties
\begin{equation*}
    V^{(1)}_{I} \cong
    \begin{cases}
    \Delta_{I} \cap H_{-1} & |I| < 2 \\
    D^{N-|I|} & |I| \geq 2
    \end{cases}
\end{equation*}
where $D^{N-|I|}$ is some polydisc in the intersection $\Delta_I \cong S^{|I|-1}\times D^{N-|I|}$ for $|I|>1$. Similarly, we may assume that the induced open regions $\{V_{st}|t \neq s\}$ restricted to $\bbP^N \setminus N_{\epsilon_1}(H_{-2})$ yields the gluing property of the induced flag $\mathfrak{L}\cap H_{-1} \cap V_s$. Inductively, we obtain an open cover $\{V_s\}$ such that the induced open cover of $H_{-j}\setminus N_{\epsilon_j}(H_{-j-1})$, denoted by $\{V^{(j)}_s\}$, satisfies the following properties: 
\begin{equation*}
    V^{(j)}_{I} \cong
    \begin{cases}
    \Delta_{I} \cap H_{-j} & |I| < j+1 \\
    D^{N-|I|} & |I| \geq j+1
    \end{cases}
\end{equation*}
for $j=0, \dots, N$. Also, the induced open regions $\{V_{st}|t \neq s\}$ restricted to $\bbP^N\setminus N_{\epsilon_j}(H_{-j-1})$ yields the gluing property of the induced flag $\mathfrak{L} \cap H_{-j} \cap V_s$. Let's write $Z_{-j}:=\pi^{-1}(H_{-j})$ and $W_{i}:=\pi^{-1}(L_i)$. 

\begin{lem}\label{lem: topological description delta fil}
    For each $q \geq 0$, we have
    \begin{equation}\label{eq: E_1 term perverse}
        H^q_{Z_{-j}-Z_{-j-1}}(Y, \C|_{W_i-W_{i-1}}) \cong \bigoplus_{|I|=j-i+1}H^{q-2r}(Y_I, Y_{I,sm};\C) 
    \end{equation}
    where $r$ is the codimension of $Z_{-j}\cap W_i$ in $W_i$. 
\end{lem}

\begin{proof}
    We first rewrite the cohomology $H^q_{Z_{-j}-Z_{-j-1}}(Y, \C|_{W_i-W_{i-1}})$ by considering the excision principle for local cohomology groups. Then we have
\begin{equation*}
\begin{aligned}
    H^q_{Z_{-j}-Z_{-j-1}}(Y, \C|_{W_i-W_{i-1}}) &\cong H^q_{Z_{-j} \cap (Y-Z_{-j-1})}(Y-Z_{-j-1}, \C|_{W_i-W_{i-1} \cap (W-Z_{-j-1})}) \\
    & \cong \left[H^q_{Z_{-j}\cap W_{i}^\circ}(W_i^\circ) \to  H^q_{Z_{-j}\cap W_{i-1}^\circ}(W_{i-1}^\circ) \right] \\
    & \cong H^{q-2r}(Z_{-j}\cap W_i^\circ, Z_{-j}\cap W_{i-1}^\circ)
\end{aligned}
\end{equation*}
where $W_\bullet^\circ:=W_\bullet-W_\bullet \cap Z_{-j-1}$. The last isomorphism comes from the tubular neighborhood theorem. By the condition (\ref{eq:tubular pair}) on the flag $\mathfrak{H}$, we have 
\begin{equation*}
    H^{q-2r}(Z_{-j}\cap W_i^\circ, Z_{-j}\cap W_{i-1}^\circ) \cong H^{q-2r}((Z_{-j}-\pi^{-1}(N_{\epsilon_j}(H_{-j-1})) \cap W_i, (Z_{-j}-\pi^{-1}(N_{\epsilon_j}(H_{-j-1}))) \cap W_{i-1})
\end{equation*}
Now, we take Mayer-Vietoris sequence with respects to $\{U_s:=\pi^{-1}(V_s)|s=0, \cdots, N\}$. Note that over $U_I$, the gluing property yields
\begin{equation*}
    U_I \cap (Z_{-j}-\pi^{-1}(N_{\epsilon_j}(H_{-j-1}))) \cong \begin{cases}
    Y_{I, sm^{(-i+j-|I|+1)}} & |I|<j-i+1 \\
    Y_I & |I| \geq j-i+1
    \end{cases}
\end{equation*}
Therefore the Mayer-Vietoris sequence is given by 
\begin{equation}\label{eq:MV for Lemma}
\begin{aligned}
    \bigoplus_{|I|=1}H^{q-2r}(Y_{I, sm^{(-i+j)}}, Y_{I, sm^{(-i+j+1)}}) &\xrightarrow{d_1} \bigoplus_{|I|=2}H^{q-2r}(Y_{I, sm^{(-i+j-1)}}, Y_{I, sm^{(-i+j)}})\\
    & \xrightarrow{d_1} \cdots \xrightarrow{d_1} \bigoplus_{|I|=j-i+1}H^{q-2r}(Y_{I}, Y_{I, sm}) \rightarrow 0.
\end{aligned}
\end{equation}
where the differential satisfies the Mayer-Vietoris sign rule (\ref{eq:MV sign rule}). Also, for $I=\{i_1, \dots, i_k\}$, the direct summand $H^{q-2r}(Y_{I, sm^{(-i+j+1-k)}}, Y_{I, sm^{(-i+j-k)}})$ of the $k$-th term can be computed by regarding the pair $(Y_{I, sm^{(-i+j+1-k)}}, Y_{I, sm^{(-i+j-k)}})$ as subspaces of $Y_{i_1}$. Note that the choice of the index ${i_1}$ does not matter because of the topological restriction we've made (Hypothesis \ref{Hypothesis}). Then it becomes $\bigoplus_J H^{q-2r}(Y_{I\cup J}, Y_{I \cup J,sm})$ for $J \subset \{0, \dots, N\}\setminus I$ with $|J|=-i+j+1-k$. In other words, the $k$-th term is given by 
\begin{equation*}
    S_k:=\bigoplus_{|I|=k}\bigoplus_{J \subset \{0, \dots, N\}\setminus I, |J|=-i+j+1-k} H^{q-2r}(Y_{I\cup J}, Y_{I \cup J,sm})
\end{equation*}
so that the sequence (\ref{eq:MV for Lemma}) becomes
\begin{equation*}
    S_1 \xrightarrow{d} S_2 \xrightarrow{d} \cdots \xrightarrow{d} S_{j-i+1} \to 0
\end{equation*}
where each $d$ is the signed sum of the isomorphisms where the signs are determined by the Mayer-Vietoris rule. In fact, it fits into the simple combinatorial sequence 
\begin{equation*}
    0 \to S_0 \xrightarrow{d} S_1 \xrightarrow{d} S_2 \xrightarrow{d} \cdots \xrightarrow{d} S_{j-i+1} \to 0
\end{equation*}
hence the conclusion follows. 
\end{proof}

Now we rewrite the $E_1$-page of the spectral sequence for the $\delta$-filtration. Recall that ${}^\delta E_1^{l,n+a}=\bigoplus_{i+j=l}H^{n+a+l}_{Z_{-j}-Z_{-j-1}}(Y,\C|_{Y_i-Y_{i-1}})$ and the differential $d_1$ is the signed sum of the connecting homomorphisms (See \ref{diag: delta E_1 page} for the precise description). Since we choose the open cover $\{U_s|s=0, \cdots, N\}$ in the proof of Lemma \ref{lem: topological description delta fil} that is independent of $i$ and $j$, the isomorphisms in Lemma \ref{lem: topological description delta fil} respects the functoriality. Then we have 
\begin{equation*}
    {}^\delta E_1^{l,n+a} \cong \bigoplus_{i+j=l}\bigoplus_{|I|=j-i+1}H^{n+a+l-2r}(Y_I, Y_{I,sm};\C) 
\end{equation*}
where $r$ is the codimension of $Z_{-j}\cap W_i$ in $W_i$. More explicitly, we have
\begin{equation*}
    \begin{tikzcd}
          E_1^{-1, n+a}\arrow{r} & E_1^{0, n+a} \arrow{r}& E_1^{1, n+a}\\
          \bigoplus_{|I|=2}H^{n+a-1}(Y_I, Y_{I,sm}) \arrow[d, phantom, "\bigoplus"] \arrow[r] \arrow[rd]& \bigoplus_{|I|=1}H^{n+a}(Y_I, Y_{I,sm}) \arrow[d, phantom, "\bigoplus"] \arrow[rd] &  &  \\
        \bigoplus_{|I|=4}H^{n+a-3}(Y_I, Y_{I,sm})\arrow[d, phantom, "\bigoplus"] \arrow[r] \arrow[rd]& \bigoplus_{|I|=3}H^{n+a-2}(Y_I, Y_{I,sm})\arrow[d, phantom, "\bigoplus"]  \arrow[r] \arrow[rd]& \bigoplus_{|I|=2}H^{n+a-1}(Y_I, Y_{I,sm}) \arrow[d, phantom, "\bigoplus"] \\
        \vdots & \bigoplus_{|I|=5}H^{n+a-4}(Y_I, Y_{I,sm}) \arrow[r] & \bigoplus_{|I|=4}H^{n+a-3}(Y_I, Y_{I,sm})  \\
    \end{tikzcd}
\end{equation*}
The horizontal (resp. anti-diagonal) differential $d_I$ (resp. $d_{II}$) is the signed sum of the relevant connecting homomorphisms $\rho^J_I$ (resp. $(\rho^J_I)^\vee$) following the Mayer-Vietoris sign rule (\ref{eq:MV sign rule}). The differential $d_1:{}^\delta E_1^{l,n+a} \to {}^\delta E_1^{l+1,n+a}$ is given by $d_1=d_I+(-1)^ld_{II}$. Therefore,
by the assumption, we have the following equivalence of the $E_1$-page of the spectral sequences
\begin{equation*}
    (\bigoplus_{p-q=a}\Gr_F^p{}^{W(M)}E_1^{l, p+q}, d_1) \cong ({}^\delta E_1^{l,n+a},d_1)
\end{equation*}
Since both spectral sequences degenerate at the $E_2$-page, we have
\begin{equation*}
    \bigoplus_{p-q=a}\Gr_F^p\Gr^{W_{lim}}_{p+q}H^{p+q+l}(X)=\bigoplus_{p-q=a}\Gr_F^p{}^{W(M)}E_2^{l, p+q} \cong {}^\delta E_2^{l,n+a}=\Gr_{\delta}^lH^{n+a+l}(Y)
\end{equation*}

The conclusion follows from Theorem \ref{thm:cataldo lef} (3):
\begin{equation*}
    {}^\delta E_2^{l,n+a}=\Gr_{\delta}^{l}H^{n+a+l}(Y)=\Gr^P_{n+a}H^{n+a+l}(Y)
\end{equation*}
This completes the proof. 

\end{proof}
\begin{proof}[Proof of Theorem \ref{thm:P=W main}-(2)]
The proof of the item (2) is almost the same. We use the Mayer-Vietoris argument to describe the spectral sequence for the flag (=perverse) filtration. To do so, we should work with the regular cohomology groups, not compactly supported ones. The ides is to apply the well-known Poincaré duality statement for the perverse filtration on $H_c^*(\tilde{Y})$:
\begin{equation*}
    \Gr^P_{n+a}H_c^{n+a+l}(\tilde{Y}) \cong (\Gr^P_{n-a}H^{n-a-l}(\tilde{Y}))^*
\end{equation*}
In terms of the $G$-flag filtration, this is isomorphic to $(\Gr_G^{-l}H^{n-a-l}(\tilde{Y}))^*$. Therefore, it is enough to show that 
\begin{equation*}
        \bigoplus_{p-q=a}\Gr_F^p\Gr^{W}_{p+q}H^{p+q+l}(X_c) \cong (\Gr_G^{-l}H^{n-a-l}(\tilde{Y}))^*
\end{equation*}
Apply the Mayer-Vietoris argument and describe the $E_1$-page of the spectral sequence for the $G$-flag filtration. We have 
\begin{equation*}
        {}^G E_1^{-l, n-a} \cong \bigoplus_{|I|=l}H^{n-a-l}(Y_I, Y_{I,sm})
\end{equation*}
and the differential $d^G_1:{}^G E_1^{-l, n-a} \to {}^G E_1^{-l+1, n-a}$ is the signed sum of the relevant connecting homomorphisms $(\rho^J_I)$ following the Mayer-Vietoris sign rule (\ref{eq:MV sign rule}). Since the spectral sequence for the $G$-filtration degenerates at $E_2$-page, we have ${}^G E_2^{-l, n-a} \cong \Gr_G^{-l}H^{n-a-l}(\tilde{Y})$. To compute the Poincaré dual $(\Gr_G^{-l}H^{n-a-l}(\tilde{Y}))^*$, we take the dual of the $E_1$-page $({}^G E_1^{-l, n-a}, d^G_1)$, denoted by $({}^G E_1^{-l, n-a^*}, (d^G_1)^*)$. By Poincaré duality (\ref{eq: PD hybrid}), this becomes 
\begin{equation*}
    ({}^G E_1^{-l, n-a})^* \cong \bigoplus_{|I|=l}H^{n+a-l}(Y_I, Y_{I,sm})
\end{equation*}
with the induced differential, the signed sum of the relevant connecting homomorphisms $(\rho^J_I)^\vee$ following the Mayer-Vietoris sign rule (\ref{eq:MV sign rule}). By the assumption (Conjecture \ref{conj:mirrorP=W}), we have the mirror equivalence of the $E_1$-page of the spectral sequences
\begin{equation*}
    (\bigoplus_{p-q=a}\Gr_F^p{}^{W}E_1^{l, p+q}, d_1) \cong (({}^GE_1^{l,n+a})^*,(d^G_1)^*)
\end{equation*}
As both spectral sequences degenerates at the $E_2$-page, we get the conclusion. 
\end{proof}

\section{Toric degeneration}\label{sec:Toric degeneration}
In this section, we propose how one can see the gluing construction in the case of semi-stable toric degenerations. 

\subsection{Backgrounds on toric varieties}
We recollect some backgrounds about toric varieties. We refer for more details to \cite{Toricvarieties}. Let $N$ and $M$ be dual lattices of rank $n$ with the natural bilinear pairing $\langle -,- \rangle:N \times M \to \Z$. We write $N_\R:=N \otimes_\Z \R$ and $M_\R=M \otimes_\Z \R$. A rational convex polyhedra cone (simply call \textit{cone}) $\alpha$ in $N_\R$ is a convex cone generated by finitely many vectors in $N$. Associated to a cone $\alpha$, one can construct an affine toric variety $X_\alpha:=\Spec(\C[\alpha^\vee \cap M])$ where $\alpha^\vee \subset M_\R$ is a dual cone of $\alpha$ given by 
    \begin{equation*}
        \alpha^\vee=\{v \in M_\R| \langle v, u\rangle \geq 0 \text{ for all } u \in \alpha\}.
    \end{equation*}
Such affine toric varieties can be glued to produce more general toric varieties. This gluing data is combinatorially encoded in a \textit{fan} $\Sigma \subset N_\R$ which is a collection of cones such that 
    \begin{enumerate}
        \item each face of a cone in $\Sigma$ is also a cone in $\Sigma$
        \item the intersection of two cones in $\Sigma$ is a face of each.
    \end{enumerate}
    Given a fan $\Sigma$, we define a toric variety $X=X_\Sigma$ by gluing the affine toric varieties $X_\alpha:=\Spec(\C[\alpha^\vee \cap M])$: Two affine toric varieties $X_\alpha$ and $X_\beta$ are glued over $X_{\alpha\beta}:=\Spec(\C[(\alpha\cap\beta)^\vee \cap M])$. We call $\{X_\alpha\}$ a toric chart of $X_\Sigma$. If $|\Sigma|=N_\R$, it is called \textit{complete}, and the corresponding toric variety $X_\Sigma$ is compact. 
    
    Let $\Sigma[1]$ be a collection of integral primitive ray generators of $\Sigma$. We consider the lattice morphism $g:N \to \Z^{\Sigma[1]}$ given by $g(v)=(\langle v, \rho\rangle)_{\rho \in \Sigma[1]}$. This induces a short exact sequence
    \begin{equation*}
        0 \to N \xrightarrow{g} \Z^{\Sigma[1]} \to A_{n-1}(X_\Sigma) \to 0
    \end{equation*}
    where $\Z^{\Sigma[1]}$ is the set of torus invariant Weil divisors and $A_{n-1}(X_\Sigma)$ is the Chow group of $X_\Sigma$. Applying the functor $\Hom(-, \C^*)$ to the above sequence, we get a short exact sequence
    \begin{equation}\label{eq:SES toric chart}
        1 \to G \to (\C^*)^{\Sigma[1]} \to M \otimes_\Z \C^* \to 1
    \end{equation}
    
    Let $\{x_\rho\}_{\rho \in \Sigma[1]}$ be a standard basis of rational functions on $\C^{\Sigma[1]}$ and $V$ be the vanishing locus of $\{\prod_{\rho \notin \sigma}x_\rho|\sigma \subset \Sigma\}$. The sequence (\ref{eq:SES toric chart}) shows that $G$ acts naturally on $\C[(x_\rho)_{\rho \in \Sigma[1]}]$ and leaves $V$ invariant. Then the toric variety $X_\Sigma$ is the quotient $(\C[(x_\rho)_{\rho \in \Sigma[1]}] \setminus V)//G$ and the homogeneous coordinate ring of $X_\Sigma$ is equipped with the grading given by the action of $G$. The sublocus of $X_\Sigma$ corresponding to $D_\rho=\{x_\rho=0\}$ is exactly the torus invariant divisors associated to the ray generator $\rho$. A torus invariant divisor $D=\sum_{\rho \in \Sigma[1]} a_\rho D_\rho$ is Cartier if and only if there is some piecewise linear function $\rho$ on $M_\R$ which linear on the cones of $\Sigma$, which takes the integral values on $M$. 
    
     A (compact) rational convex polytope $\Delta$ in $M_\R$ is the convex hull of finite number of points. We say $\Delta$ is a \textit{lattice polytope} if every vertex of $\Delta$ is in $M$. For example, a lattice polytope is given by
    \begin{equation*}
        \Delta=\{v \in M_\R|\langle v, n_i \rangle \geq -a_i, n_i \in N, a_i \in \Z \text{ for }i=1,\dots, s\}
    \end{equation*}
    A \textit{$l$-face} $\sigma$ is the intersection of $\Delta$ with $n-l$ supporting hyperplanes, and we will it denote by $\sigma \prec \Delta$.  We also write $\Delta[l]$ for the collection of $l$-faces of $\Delta$. In particular, a $0$-face, a $1$-face and a $(n-1)$-face are called a vertex, an edge and a facet of $\Delta$, respectively. For each cone $\sigma \prec \Delta$, the cone $\alpha_\sigma$ dual to $\sigma$ is defined by 
    \begin{equation*}
        \alpha_\sigma=\{u \in N_\R|\langle v, u\rangle \leq \langle v', u \rangle \text{ for all } v \in \sigma \text{ and } v' \in \Delta\}
    \end{equation*}
    A collection of dual cones $\alpha_\sigma$ forms a fan $\Sigma^{\Delta}$, called a normal fan of $\Delta$, and we write $X_\Delta$ for the associated toric variety.
    
    To a lattice polytope $\Delta$, the associated toric variety $X_\Delta$ comes with the divisor $D_\Delta=-\sum_{\rho}a_\rho D_\rho$ (or simply denoted by $D$) where the sum is taken over all facets $\rho \prec \Delta$. Equivalently, we get a support function of $D$, a piecewise linear function $\phi_D$ such that $\phi_D(v_\rho)=-a_\rho$ for the privimite vector $v_\rho$ to the dual cone of the face $\rho \prec \Delta$.
    Let $\Delta_D=\{u \in M_\R|u \geq \phi_D \text{ on } N_\R\}$. Geometrically, $\Delta_D \cap M$ generates the space of sections of the line bundle $\cO_X(D)$. Note that $D$ is trivial, generated by sections and ample if and only if $\phi_D$ is affine, convex, and strictly convex, respectively.
    
    A polytope $\Delta \subset M_\R$ is called \textit{simlicial}, if there are exactly $n$ edges at each vertex and the primary vectors at each vertex span $M_\R$ as a vector space. A fan $\Sigma$ in $N_\R$ is simplicial if all the maximal cones in $\Sigma$ is simplicial. In particular, if the primary vectors span the lattice, then it is called \textit{non-singular}.
     \begin{prop}
    If $\Delta$ is simplicial (resp. non-singular), then $X_\Delta$ is an orbifold (resp. manifold).
    \end{prop}
    \subsection{Batyrev-Borisov mirror and Givental's LG model}\label{sec:BBmirror}
    We introduce Batyrev-Borisov mirror pairs \cite{BatyrevBorisovmirror}. Let $\Delta$ be a simplicial lattice polytope in $M_\R$. A polar dual of the polytope $\Delta$ is defined to be $\Delta^\circ:=\{u \in N_\R|\langle u, v \rangle \geq -1 \text{ for }v \in \Delta\}$. A lattice polytope is called \textit{reflexive} if its polar dual $\Delta^\circ$ is also a lattice polytope. The equivalent condition is that $\Delta$ has only one interior lattice point, denoted by $0_M$. Geometrically, the associated toric variety $X_\Delta$ is a Gorenstein Fano variety. From know on, we fix a reflexive polytope $\Delta \subset M_\R$ and write $\Sigma_\Delta$ for the fan over the facets of $\Delta$ and $\bbP_\Delta:=X_{\Sigma_\Delta}$ for the associated toric variety. Note that $\Sigma_\Delta$ is also the normal fan of the polar dual $\Delta^\circ$, so we have $\bbP_\Delta=X_{\Sigma_\Delta}=X_{\Sigma^{\Delta^\circ}}=X_{\Delta^\circ}$.
    
    Consider a general Calabi-Yau hypersurface $V_\Delta$ of $\bbP_\Delta$. Since $\bbP_\Delta$ is an orbifold, the hypersurface $V_\Delta$ may have singularities. We assume that the hypersurface $V_\Delta$ is $\Delta$-regular, meaning that the singular locus of $V_\Delta$ is induced from the singular locus of ambient space $\bbP_\Delta$. Then one may desingularize $V_\Delta$ by taking a partial resolution of $\bbP_\Delta$. Such resolution is given by a refinement $\widetilde{\Sigma}_\Delta$ of the fan $\Sigma_\Delta$ whose cone is contained in a cone of $\Sigma_\Delta$. In this case, to $\Sigma_\Delta$, one can add all rays pointing to the elements in $\partial \Delta \cap M$ to obtain $\widetilde{\Sigma}_\Delta$. Batyrev shows that the induced resolution $f: X_{\widetilde{\Sigma}_\Delta} \to \bbP_\Delta$ is crepant and this is called a maximal projective crepant partial (MPCP) resolution of $\bbP_\Delta$\cite[Section 2.2]{Batyrevdualpolyhedra}. In particular, if $X:=f^*(V_\Delta)$ is smooth, then we say $\Delta$ satisfies MPCS (maximal projective crepant smooth) resolution condition. Note that this condition always holds for $n \leq 4$ \cite[Section 2.2]{Batyrevdualpolyhedra}. Similarly, consider the dual construction for $\Delta^\circ$ and write $X^\vee$ for a MPCP resolution of $V_{\Delta^\circ}$. 
    
    \begin{thm}\label{thm:batyrevBorisovhypersurface}\cite{Batyrevdualpolyhedra}\cite{BatyrevBorisovmirror}
        The pair of Calabi-Yau hypersurfaces $(V_\Delta, V_{\Delta^\circ})$ (or $(X, X^\vee)$) satisfies (stringy) Hodge number mirror relation. We call this pair a Batyrev-Borisov mirror pair. 
    \end{thm}
    Theorem \ref{thm:batyrevBorisovhypersurface} is generalized to general Calabi-Yau complete intersections. 
    
    \begin{defin}
        A nef partition of $\Delta$ is a partition of $\Delta[0]$ into subsets $E_1, \dots, E_{k+1}$ such that there exists integral convex $\Sigma_\Delta$-functions on $M_\R$ such that $\phi_i(v)=1$ for $v \in E_i$ and $\phi_i(v)=0$ otherwise. 
    \end{defin}

    Let $\Delta_i=\mathrm{Conv}(0_M \cup E_i)$ be the convex hull of $E_i$. We have the Minkowski decomposition $\Delta=\Delta_1 +\cdots+ \Delta_{k+1}$ and $\cap_{i=1}^{k+1} \Delta_i=\{0_M\}$. Note that the nef partition allows one to decompose the anti-canonical divisor into $k+1$-line bundles
    \begin{equation*}
        \cO_{\bbP_\Delta}(\sum_{\rho \in \Delta[0]}D_\rho)=\cO_{\bbP_\Delta}(\sum_{\rho \in E_1}D_\rho)\otimes\cdots \otimes  \cO_{\bbP_\Delta}(\sum_{\rho \in E_{k+1}}D_\rho)
    \end{equation*}
    This induces a general Calabi-Yau complete intersection of codimension $k+1$, denoted by $V_\Delta$ as before. If we take the MPCP resolution $f:X_{\widetilde{\Sigma}_\Delta} \to \bbP_\Delta$, we have $\cO_{X_{\widetilde{\Sigma}_\Delta}}(\sum_{\rho \in \partial \Delta_i\setminus \{0\}}D_\rho)=f^*\cO_{\bbP_\Delta}(\sum_{\rho \in E_i}D_\rho)$ for each $i$. If we put $\Delta$-regularity conditions on each hypersurfaces, the MPCP resolution $f$ induces a (partial) desingularization of $V_\Delta \subset \bbP_\Delta$, denoted by $X$. 
    
    To describe mirror dual, we define the polar dual of $\Delta_i$ to be 
    \begin{equation*}
        \nabla_i:=\{u \in N_\R| \langle u, v \rangle \geq -\phi_i(v)\}
    \end{equation*}
    By construction, we have $\Delta^\circ=\nabla_1 + \cdots + \nabla_{k+1}$. We denote the convex hull of $\nabla_i$'s by $\nabla=\mathrm{Conv}(\nabla_1 \cup \dots \cup \nabla_{k+1})$. In general, $\nabla$ is strictly contained in $\Delta^\circ$. We take a general Calabi-Yau complete intersection $V_\nabla$ in $\bbP_\nabla$ (instead of $\bbP_{\Delta^\circ}$) and its (partial) desingularization $X^\vee$. 
    
    \begin{thm}\cite{BatyrevBorisovmirror}
        The pair of general Calabi-Yau complete intersections $(V_\Delta, V_{\nabla})$ (or $(X, X^\vee)$) satisfies (stringy) Hodge number mirror relation. We call this pair a Batyrev-Borisov mirror pair. 
    \end{thm}
    
    Next, consider a general complete intersection of the first $k$ divisors in $\bbP_\Delta$, denoted by $V_{(k)}$. This is a quasi-Fano variety with an anti-canonical divisor $V_\Delta$. In \cite{Giventalmirrorcomplete}, A.Givental introduces a mirror construction of the pair $(V_{(k)}, V_\Delta)$. Let $(\C^*)^n=\C[x_1^{\pm}, \dots, x_n^{\pm}]$ be a $n$-dimensional affine torus. The Givental's LG model is given by a pair $(Y_0, \sw_0)$ where 
     \begin{equation*}
         \begin{aligned}
            &Y_0=\bigcap_{i=1}^k \{ \sum_{\rho \in \Delta_i \cap M } a_\rho x_\rho=0\} \\
            &\sw_0=\sum_{\rho \in \Delta_{k+1} \cap M}a_\rho x_\rho
         \end{aligned}
     \end{equation*}
     where $a_\rho$ are general complex coefficients and $x_\rho=\prod_{i=1}^n x_i^{\langle e_i,\rho\rangle}$ where $\{e_1, \dots, e_n\}$ is the basis of $N$. One can get a fiberwise compactification of the LG model $(Y_0, \sw_0)$ by embedding each fiber into $X_{\widetilde{\Sigma}_\nabla}$ and view $\sw_0$ as a pencil of hypersurfaces. The reason is that Laurant polynomial $\sum_{\rho \in \Delta_i \cap M}a_\rho x_\rho$ can be polarized to a section of $\cO_{X_{\widetilde{\Sigma}_\nabla}}(\sum_{\sigma \in \nabla_i \cap N}D_\sigma)$ that corresponds to an integral point of $M$. Then a fiber at $\lambda \in \C$ is compactified to the following:
     \begin{equation}\label{eq:LGcompact}
         \begin{aligned}
             \sum_{\rho \in \Delta_i \cap M} a_\rho \prod_{\sigma \in \nabla \cap N}z_{\sigma}^{\langle \sigma, \rho \rangle-\sigma^i_{min}}=0\\
             \lambda \prod_{\sigma \in \nabla_{k+1} \cap N\setminus 0_N}z_\sigma-\sum_{\rho \in \Delta_{k+1}\cap M\setminus 0_M} a_\rho \prod_{\sigma \in \nabla \cap N}z_{\sigma}^{\langle \sigma, \rho \rangle-\sigma^{k+1}_{min}}=0
         \end{aligned}
     \end{equation}
     where $z_\sigma$ are homogenous coordinate of $X_{\widetilde{\Sigma}_\nabla}$ and $\sigma^i_{\min}=\min\{\langle \sigma, \rho \rangle | \rho \in \Delta_i\}$ for all $i$. 
     
     \begin{thm}\cite[Theorem 3.2.6]{HarderThesis}\label{thm:Harder compactification LG}
         Let $V_{(k)}$ be a complete intersection quasi-Fano variety associated to a $(k+1)$-partite nef partition of a polytope $\Delta$. Then the partial compactification of $(Y_0, \sw_0)$ to the total space of the pencil $\mathcal{P}(\lambda) \subset X_{\widetilde{\Sigma}_\nabla} \times \bbA^1$ has at worst terminal Gorenstein singularities, is log Calabi-Yau and admits a compactification $Z$ where $\sw$ extends to a map $f:Z \to \bbP^1$ so that $f^{-1}(\infty)$ is a normal crossing union of varieties with at worst toroidal singularities. 
     \end{thm}
     
     In Theorem \ref{thm:Harder compactification LG}, the total space of pencil $\mathcal{P}(\lambda)$ is given by the same equations in (\ref{eq:LGcompact}) except that $\lambda$ is shifted by $\lambda+a$ for some non-zero constant $a$. 
     
    Suppose we have a partition of $E_{k+1}$ into subsets $F_1, \dots, F_r$. In case that this partition is nef, then $\Delta$ admits the $(k+r)$-partite nef partition a general complete intersection quasi-Fano variety $V_{(k)}$ comes with normal crossing anti-canonical divisor. We can generalize the Givental's mirror construction of $V$ with a simple normal crossing anti-canonical divisor whose component corresponds to a generic section of $\cO_{\bbP_{\Delta}}(\sum_{\rho \in \Delta_{k+j} \cap M \setminus \{0_M\}} D_\rho)$ for each $j=1, \dots, r$. Then mirror LG model is given by a pair $(Y_0, h=(h_1, \dots, h_r):Y_0 \to \C^r)$ where 
   \begin{equation}\label{eq:givental mirror hybrid}
         \begin{aligned}
            &Y_0=\bigcap_{i=1}^k \{ \sum_{\rho \in \Delta_i \cap M } a_\rho x_\rho=0\} \\
            &h=(h_1, \cdots, h_r):Y_0 \to \C^r, \quad h_j=\sum_{\rho \in \Delta_{k+j} \cap M}a_\rho x_\rho
         \end{aligned}
     \end{equation}
    Similar to the $r=1$ case, one can get a family of Calabi-Yau varieties over $(\bbP^1)^r$. In other words, a fiber at $\lambda=(\lambda_1, \dots, \lambda_r)$ is compactified to the subvariety in $X_{\widetilde{\Sigma}_\nabla}$ cut out by
     \begin{equation*}
         \begin{aligned}
             \sum_{\rho \in \Delta_i \cap M} a_\rho \prod_{\sigma \in \nabla \cap N}z_{\sigma}^{\langle \sigma, \rho \rangle-\sigma^i_{min}}=0 & \qquad 1 \leq i \leq k \\
             \lambda_j \prod_{\sigma \in \nabla_{k+j} \cap N \setminus 0_N}z_\sigma-\sum_{\rho \in \Delta_{k+j}\cap M \setminus 0_M} a_\rho \prod_{\sigma \in \nabla \cap N}z_{\sigma}^{\langle \sigma, \rho \rangle-\sigma^{k+j}_{min}}=0 & \qquad 1 \leq j \leq r
         \end{aligned}
     \end{equation*}
     where $z_\sigma$ are homogenous coordinate of $X_{\widetilde{\Sigma}_\nabla}$ and $\sigma^i_{\min}=\min\{\langle \sigma, \rho \rangle | \rho \in \Delta_i\}$ for all $i=1, \dots, k+r$. 

    However, we will encounter the case that the partition of $E_{k+1}$ is not nef. One naive way to obtain a hybrid LG potential is to consider $\lambda$ as the sum of $r$-variables $\lambda_1, \dots, \lambda_r$ and consider the subvariety in $X_{\widetilde{\Sigma}_\nabla} \times \C^r$ where 
      \begin{equation*}
         \begin{aligned}
             \sum_{\rho \in \Delta_i \cap M} a_\rho \prod_{\sigma \in \nabla \cap N}z_{\sigma}^{\langle \sigma, \rho \rangle-\sigma^i_{min}}=0 & \qquad 1 \leq i \leq k \\
             \lambda_j \prod_{\sigma \in \nabla_{k+i} \cap N \setminus 0_N}z_\sigma-\sum_{\rho \in F_j} a_\rho \prod_{\sigma \in \nabla \cap N}z_{\sigma}^{\langle \sigma, \rho \rangle-\sigma^{k+1}_{min}}=0 & \qquad 1 \leq j \leq r
         \end{aligned}
     \end{equation*}
    Now the question whether this family has nice enough singularities and provides a mirror hybrid LG model of $V$ with a simple normal crossing anti-canonical divisor.

\subsection{Semi-stable partition}\label{sec:semistable partition}
We start with reviewing a semi-stable toric degeneration introduced in \cite{Semistabledegenerationtoric}. Fix a $n$-dimensional simplicial polytope $\Delta$. Let $\Gamma$ be a partition of the polytope $\Delta$ into smaller polytopes $\{\Delta_{(i)}\}$. We say the partition $\Gamma$ is \textit{simplicial} if all subpolytopes $\{\Delta_{(i)}\}$ are simplicial polytopes. We define $\sigma$ to be $l$-face of $\Gamma$, denoted by $\sigma \prec \Gamma$, if $\sigma$ is a $l$-face of $\Delta_{(i)}$ for some $i$. 

\begin{defin}\label{def:semistable parition}
A simplicial partition $\Gamma$ is \textit{semi-stable} if the following conditions hold:
\begin{enumerate}
    \item each vertex of $\Delta$ belongs to only one of $\Delta_{(i)}$'s;
    \item for any $l$-face $\sigma \prec \Gamma$ and $k$-face $\tau \prec \Delta$ with $\sigma \subset \tau$, then there are exactly $k-l+1$ $\Delta_{(i)}$'s such that $\sigma \prec \Delta_{(i)}$.
\end{enumerate} 
\end{defin}

From now on, we distinguish vertices in $\Gamma$ from those of $\Delta$: when we say $p$ is a vertex of $\Gamma$, it means that $p$ is the one that is not in $\Delta$. The restriction of $\Gamma$ to some face $\sigma \prec\Delta$ is defined to be the partition induced by $\{\Delta_{(i)} \cap \sigma\}$, and denote it by $\Gamma \cap \sigma$. 

By definition, if $\cap_{k=1}^l\Delta_{(i_k)} \neq \emptyset$, then it has dimension $n-l+1$. Also, for any vertex $p$ of $\Gamma$, there are exactly $n+1$-edges $\sigma_0, \dots, \sigma_n$ of $\Gamma$ such that $p \prec \sigma_i$. This allows us to define a dual simplicial complex $K_\Gamma$ whose vertex set is the set of polytopes $\Delta_{(j)}$'s in the partition $\Gamma$. For example, $K_\Gamma$ is $l$-simplex if and only if there is a $l$-face of $\Delta$ that contains all vertices of $\Gamma$. 

\begin{defin}
    A vertex $p \prec \Gamma$ is non-singular if $p$ is non-singular in one (thus all) sub-polytope containing $p$. A semi-stable partition $\Gamma$ is non-singular if all of vertices are non-singular. 
\end{defin}

    \begin{defin}
        A lifting of $\Delta$ by a semi-stable partition $\Gamma$ is a triple $(\Tilde{\Delta}, \Tilde{M}, \pi)$ where $\Tilde{\Delta}$ is a lattice polytope of $\Tilde{M}$ and $\pi:\tilde{M} \to M$ is a surjective morphism satisfying the following condition: for $\tilde{\sigma} \prec \tilde{\Delta}$, either $\pi_*(\tilde{\sigma}) \prec \Delta$ or $\pi_*(\tilde{\sigma}) \prec \Gamma$ where $\pi_*:\tilde{M}_\R \to M_\R$ is the induced map from $\pi$. If $\pi_*(\tilde{\sigma}) \prec \Gamma$, $\pi_*(\tilde{\sigma})$ is said to be a lift of $\pi_*(\tilde{\sigma})$. Also, the lifting is called non-singular if all polytopes involved are non-singular. 
    \end{defin}
    
    \begin{prop}\cite[Proposition 3.12]{Semistabledegenerationtoric}
        For a non-singular semi-stable partition $\Gamma$, there exists a concave integral piecewise linear function $F_\Gamma$ on $\Delta$ that is linear on each face of $\Gamma$. 
    \end{prop}
    
    One can take minimal integral lifting of $F_\Gamma$ and we denote it by $F$. 

    \begin{thm}\cite[Theorem 3.13]{Semistabledegenerationtoric}
        Let $\pi$ be the projection $\Z \oplus M \to M$ and $\Tilde{\Delta}=\{(y,x)|y \geq F(x)\} \subset \R \times \Delta \subset \R \times M_\R$. Then $(\Tilde{\Delta}, \Z \oplus M, \pi)$ is a lifting of $\Delta$ by $\Gamma$. If $\gamma$ is a non-singular partition of a non-singular polytope, then the lifting is non-singular.
    \end{thm}
    \begin{thm}\cite[Theorem 4.1]{Semistabledegenerationtoric}
        Suppose $\Delta$ and $\Gamma$ are both non-singular. Then there exists a semi-stable degeneration $p:X_{\tilde{\Delta}}\to \C$ of $X_\Delta$ to $p^{-1}(0)$. The dual graph $G$ of the degeneration fiber $p^{-1}(0)$ is isomorphic to $K_\Gamma$, and each component in $p^{-1}(0)$ is the toric variety defined by the corresponding subpolytope in $\Gamma$.
    \end{thm}
    
    The semi-stable degeneration $p:X_{\tilde{\Delta}}\to \C$ induces a semi-stable degeneration of a Calabi-Yau hypersurface $X$ in $X_\Delta$ whose degeneration fiber consists of a generic hypersurface $X_i$ of $X_{\Delta_{(i)}}$ in the linear system $|D_{\Delta_{(i)}}|$. In other words, the normal crossing variety $X_c=\cup_{i=0}^l X_i$ is $d$-semistable of type $l+1$ and $X$ is its smoothing. 
    
    Next, we assume that the non-singular polytope $\Delta \subset M_\R$ is reflexive and write the unique lattice point by $0_M$ and call it the origin.
    \begin{defin}
        A semi-stable partition $\Gamma$ is \textit{central} if $0_M \prec \Delta_{(i)}$ for all $i$. 
    \end{defin}

    Let $\Gamma$ be a non-singular, central, semi-stable partition. If $K_\Gamma$ is $l$-dimensional, there is a unique codimension $l$ linear subspace $L \subset M_\R$ passing through the origin such that $\Delta \cap L \prec \Delta_{(i)}$ for all $i$ as a $n-l$-face. Consider $l$ primitive vectors given by the intersection of each $\Delta_{(ik)}:=\Delta_{(i)} \cap \Delta_{(k)}$ with the orthogonal complement $L^\perp$ of $L$. Note that the restriction $\Delta \cap L^\perp$ is also reflexive and simplicial while the induced partition is not necessarily semi-stable. Thus, there are $l+1$ such primitive vectors $v_0, \dots, v_l$ such that for each $i$, $\Delta_{(i)}$ contains all $v_j$ except $j=i$. 

    Next, consider the dual reflexive polytope $\Delta^\circ \subset N_\R$ whose dual fan is $\Sigma_{\Delta} \subset M_\R$. For simplicity, we assume that $\Delta^\circ$ satisfies MPCS resolution condition. We explain what the central semi-stable partition corresponds to in the dual picture. Let $v_0, \dots, v_l$ be primitive vectors introduced above. We define new fans $\Sigma'_\Gamma \subset \Sigma_\Gamma \subset M_\R$ where 
    \begin{enumerate}
        \item $\Sigma'_\Gamma$ is generated by primitive vectors in all $\Delta_{(ik)}$'s and all $v_i$'s;
        \item $\Sigma_\Gamma$ is a subfan of $\Sigma'_\Gamma$ that is generated by all primitive vectors lying in $\Delta \cap L$ and all $v_i's$
    \end{enumerate}
    Furthermore, we consider the fan $\Sigma':=\Sigma_\Delta \cup \Sigma'_\Gamma$. Geometrically, this refinement amounts to taking a (maximal) projective crepant partial resolution of $X_{\Sigma_{\Delta}}$, denoted by  $\phi_\Delta:  X_{\Sigma'} \to X_{\Sigma_\Delta}$. We also have a blow-down map $\phi_\Gamma:X_{\Sigma'} \to X_{\Sigma_\Gamma}$ which is not necessarily crepant. Consider the projection of the lattices $\Pi_v:M \to M_v$ where $M_v$ is the sublattice of $M$ generated by $v_i$'s. This provides a surjective toric morphism $\pi_v:X_{\Sigma_\Gamma} \to X_{\Sigma_v}$ where $\Sigma_v \subset M_{v,\R}$ is the fan generated by all $v_i$'s. In fact, $\pi_v$ is a trivial fibration whose fiber is a toric variety associated to a fan generated by all primitive vectors in $\Delta \cap L$, denoted by $X_L$. In summary, we have the following diagrams of toric morphisms
    \begin{equation*}
        \begin{tikzcd}
        & X_{\Sigma'} \arrow[ld, swap, "\phi_\Delta"] \arrow[rd,  "\phi_\Gamma"] \arrow[ bend left=20, rrd, "\pi_\Gamma"] && \\
        X_{\Sigma_\Delta} && X_{\Sigma_\Gamma} \arrow[r, "\pi_v"] & X_{\Sigma_v}
        \end{tikzcd}
    \end{equation*}
    where $\pi_\Gamma:=\phi_\Gamma \circ \pi_v$. Since $\pi_\Gamma:X_{\Sigma'} \to X_{\Sigma_v}$ is the composition of toric morphisms, a generic fiber is the toric variety $X_L$. In fact, over the open dense torus $(\C^*)^l \subset X_{\Sigma_v}$, the morphism $\pi_\Gamma$ is trivial fibration with fiber $X_L$. We also present the coordinate description. Consider the homogeneous coordinates of $X_{\Sigma'}$ (resp. $X_{\Sigma_v}$), $(z_\sigma|\sigma \in \Sigma')$ (resp. $(z_{\sigma_0}, \dots, z_{\sigma_l})$). In terms of these coordinates, $\pi_\Gamma$ is given by 
    
    \begin{equation*}
\pi_\Gamma=\left(\prod_{\Pi_v(\sigma)=\sigma_0} z_\sigma, \dots, \prod_{\Pi_v(\sigma)=\sigma_l}z_\sigma\right)
    \end{equation*}
    Let's take a generic anti-canonical divisor $Y$ in $X_{\Sigma'}$ and the induced fibration $\pi:=\pi_\Gamma|_Y$. Recall that we assume $\Delta^\circ$ satisfies MPCS resolution condition so that $Y$ is non-singular. Also, it is clear that a generic fiber of $\pi$ is Calabi-Yau as this is a non-singular general member of the anti-canonical divisor of $X_L$. 

    Since $\Gamma$ is non-singular, the base $X_{\Sigma_v}$ is isomorphic to the projective space $\bbP^l$ and we write $\{z_i:=z_{v_i}|i=0, \dots, l\}$ for the homogeneous coordinates. For each $i$, we take a polydisc $\Delta_i=\{|z_j| \leq |z_i| |j \neq i \} \subset \bbP^l$ and set $Y_i:=\pi^{-1}(\Delta_i)$ and $h_i:=\pi|_{Y_i}:Y_i \to \Delta_i$. It follows from the genericity condition on $Y$ that the restriction of $\pi:Y \to \bbP^l$ over each boundary component of $\Delta_i$ is locally trivial with smooth fibers and intersects transversally to each other. This implies that $(Y_i, h_i:Y \to \Delta_i)$ is a hybrid LG model.

    Recall that on the degeneration side, the semi-stable partition $\Gamma$ provides a semi-stable degeneration of a non-singular Calabi-Yau hypersurface of $X_\Delta$. Each irreducible component $X_i$ of the degeneration fiber $X_c=\cup_{i=0}^l X_i$ is a general hypersurface of the toric variety $X_{\Delta_{(i)}}$ that is determined by all facets in the facets of $\Delta$. 

    \begin{conj}\label{Conj:hybrid LG model}
        For each $i$, the hybrid LG model $(Y_i,h_i: Y_i \to \Delta_i)$ is mirror to the pair $(X_i,\cup_{j \neq i}X_{ij})$. 
    \end{conj}
    
\begin{rem}
    One may apply the same construction without imposing MPCS resolution condition on $\Delta^\circ$. In this case, $Y$ becomes an orbifold one needs more general notion of hybrid LG models. 
\end{rem}

    \begin{example}
        Consider the reflexive square $\Delta$ with the semistable partition $\Gamma$ given by the vertical line
        
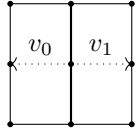
\begin{equation*}
             \begin{tikzpicture}[scale=0.6]
\foreach \x in {0,1,2}
   \foreach \y in {0,1,2} 
      \draw[fill] (4/3*\x,4/3*\y) circle (1.5pt) coordinate (m-\x-\y);

\draw (m-0-0) -- (m-2-0) -- (m-2-2) -- (m-0-2) -- cycle;
\draw[thick] (m-1-0) -- (m-1-2);
\draw[dotted,->] (m-1-1) -- (m-0-1) node [midway, above] {$v_0$};
\draw[dotted,->] (m-1-1) -- (m-2-1) node [midway, above]{$v_1$};

\end{tikzpicture}
        \end{equation*}
    where dotted arrows are primitive vectors $v_0$ and $v_1$. This describes a semi-stable degeneration of a non-singular Calabi-Yau hypersurface of $X_\Delta=\bbP^1\times \bbP^1$, which is elliptic curve, into the union of two rational curves $X_0$ and $X_1$ intersecting over two points $X_{01}$. Consider the fans $\Sigma_\Delta$, $\Sigma'$ and $ \Sigma_v$ described below (Figure 2):

\begin{center}
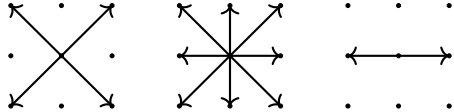

 \begin{tikzpicture}[scale=0.50]
\foreach \x in {0,1,2}
   \foreach \y in {0,1,2} 
      \draw[fill] (4/3*\x,4/3*\y) circle (1.5pt) coordinate (m-\x-\y);
\draw[thick,->] (m-1-1) -- (m-0-2);
\draw[thick,->] (m-1-1) -- (m-2-2);
\draw[thick,->] (m-1-1) -- (m-2-0);
\draw[thick,->] (m-1-1) -- (m-0-0);
\end{tikzpicture}
\qquad
\begin{tikzpicture}[scale=0.50]
\foreach \x in {0,1,2}
   \foreach \y in {0,1,2} 
      \draw[fill] (4/3*\x,4/3*\y) circle (1.5pt) coordinate (m-\x-\y);
\draw[thick,->] (m-1-1) -- (m-0-2);
\draw[thick,->] (m-1-1) -- (m-1-2);
\draw[thick,->] (m-1-1) -- (m-2-1);
\draw[thick,->] (m-1-1) -- (m-1-0);
\draw[thick,->] (m-1-1) -- (m-0-1);
\draw[thick,->] (m-1-1) -- (m-2-2);
\draw[thick,->] (m-1-1) -- (m-2-0);
\draw[thick,->] (m-1-1) -- (m-0-0);
\end{tikzpicture}  \qquad
\begin{tikzpicture}[scale=0.50]
\foreach \x in {0,1,2}
   \foreach \y in {0,1,2} 
      \draw[fill] (4/3*\x,4/3*\y) circle (1.5pt) coordinate (m-\x-\y);
\draw[thick,->] (m-1-1) -- (m-0-1);
\draw[thick,->] (m-1-1) -- (m-2-1);
\end{tikzpicture}
\captionof{figure}{Fans $\Sigma_\Delta$, $\Sigma'$ and $ \Sigma_v$ from the left}
\end{center}
     Geometrically, $X_{\Sigma'}$ is a blow up of $\bbP^1 \times \bbP^1$ along the 4 corners and $X_{\Sigma_v} \cong \bbP^1$. The morphism $\pi_\Gamma:X_{\Sigma'} \to X_{\Sigma_v}$ is the one given by the projection to the first factor. In terms of the homogenous coordinate $\{z_\sigma|\sigma \in \Sigma'[1]\}$, this is given by 
    \begin{equation}\label{eq:homo coor p1}
\pi_\Gamma=\left[z_{\sigma(-1,1)}z_{\sigma(-1,0)}z_{\sigma(-1,-1)}: z_{\sigma(1,1)}z_{\sigma(1,0)}z_{\sigma(1,1)}\right].
    \end{equation}
    
    Take a generic section $Y$ of $|-K_{X_{\Sigma'}}|$ that is not singular over the locus $\{|z_0|=|z_1|\} \in \bbP^1 \cong X_v$. Then the induced fibration $\pi:Y \to \bbP^1$ becomes a double cover with 4 ramification points. Note that two of them lie near $0 \in \bbP^1$ while the other two points lie near $\infty \in \bbP^1$. For $i=0,1$, we take $Y_i:=\pi^{-1}(\Delta_i)$ and $\sw_i:=\pi|_{Y_i}:Y_i \to \Delta_i$. We show that the pair $(Y_i, \sw_i)$ is mirror to the pair $(X_i, X_{01})$. To describe mirror of $X_0$ (the parallel argument works for $X_1$), we first make the rectangle $\Delta{(0)}$ reflexive by shifting the middle vertical line (a facet of $\Delta_{(0)}$) to the right by length $1$. We still denote it by $\Delta_{(0)}$. Then we apply Givental's construction to get a mirror of $X_0$. We consider the polar dual of $\Delta_{(0)}$, denoted by $\Delta^{(0)}$, and regard it belongs to $M$ to match the notation we have used. Note that we have a nef partition of $\Delta^{(0)}=\Delta^{(0)}_1 + \Delta^{(0)}_2$ where $\Delta^{(0)}_1[0]$ are the red vertices and $\Delta^{(0)}_2[0]$ is blue one in the following picture (left).
    \\
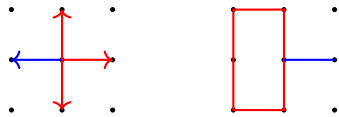
    \begin{center}
    \begin{tikzpicture}[scale=0.50]
\foreach \x in {0,1,2}
   \foreach \y in {0,1,2} 
      \draw[fill] (4/3*\x,4/3*\y) circle (1.5pt) coordinate (m-\x-\y);

% using named coordinates
\draw[thick,blue,->] (m-1-1) -- (m-0-1);
\draw[thick,red,->] (m-1-1) -- (m-2-1);
\draw[thick,red,->] (m-1-1) -- (m-1-2);
\draw[thick,red,->] (m-1-1) -- (m-1-0);
\end{tikzpicture} 
\qquad \qquad 
      \begin{tikzpicture}[scale=0.50]
\foreach \x in {0,1,2}
   \foreach \y in {0,1,2} 
      \draw[fill] (4/3*\x,4/3*\y) circle (1.5pt) coordinate (m-\x-\y);

% using named coordinates
\draw[thick, red] (m-0-0) -- (m-1-0) -- (m-1-2) -- (m-0-2) -- cycle;
\draw[thick, blue] (m-1-1) -- (m-2-1);
\end{tikzpicture}
    \end{center}
   
The right picture describes the dual partition $\nabla^{(0)}_1$ and $\nabla^{(0)}_2$. Then the LG model for $X_0$ is given by 
\begin{equation*}
    Y_0=\{\sum_{\rho \in \Delta^{(1)}_1 \cap M}a_\rho x_\rho =0 \} \subset (\C^*)^2, \qquad \sw_0=\sum_{\rho \in \Delta^{(1)}_2 \cap M}a_\rho x_\rho
\end{equation*}
 where $a_\rho$'s are general complex coefficients and $x_\rho=\prod_{i=1}^nx_i^{\langle e_i, \rho\rangle}$ where $\{e_1, \dots, e_n\}$ is the basis of $N$. A fiber is of $\sw_0$ is compactifed to 
\begin{equation*}
    \begin{aligned}
    &a_{\rho(0,0)}z_{\sigma(0,1)}z_{\sigma(-1,1)}z_{\sigma(-1,0)}z_{\sigma(-1,-1)}z_{\sigma(0,-1)}\\
    &+a_{\rho(1,0)}z_{\sigma(0,1)}z_{\sigma(0,-1)}z_{\sigma(1,0)}\\
    &+a_{\rho(0,1)}z^2_{\sigma(0,1)}z^2_{\sigma(-1,1)}z_{\sigma(-1,0)}\\
    &+a_{\rho(0,-1)}z^2_{\sigma(0,-1)}z^2_{\sigma(-1,-1)}z_{\sigma(-1,0)}=0 
    \end{aligned}
\end{equation*}
\begin{equation*}
    \lambda  z_{\sigma(1,0)}-a_{\rho(-1,0)}z_{\sigma(-1,1)}z_{\sigma(-1,0)}z_{\sigma(-1,-1)}=0
\end{equation*}
where $z_\sigma$ is the homogenous coordinate of $X_{\nabla^{(0)}}$. For generically chosen $a_\rho$'s, we see that near $\lambda=0$, this is a double cover of $\C$ with two ramification points. In fact, locally, this morphism is the same with one described in (\ref{eq:homo coor p1}). As we have the parallel argument for $\Delta_{(1)}$, we get the conclusion. 
    \end{example}
    
    There are two main issues to verify this conjecture in general. 
    \begin{enumerate}
        \item On the degeneration side, as introduced in Section \ref{sec:BBmirror}, one possible way to obtain a mirror of each pair $(X_i,\cup_{j \neq i}X_{ij})$ is to make each subpolytope $\Delta_{(i)}$ reflexive and apply Givental's construction. For each $\Delta_{(i)} \subset M_\R$, let $\Delta_{(i1)}, \dots, \Delta_{(il)}$ be the facets that is not contained in the facets of $\Delta$. If we shift the defining $l$ linear hyperplanes of these facets by length 1, then we obtain a reflexive polytope, still denoted by $\Delta_{(i)}$. We label other remaining facets by $\{\Delta_{(ij)}|j=l+1, \dots, r_i\}$ for some $r_i >0$. If we denote $\rho_\bullet$ the primitive normal vector for $\Delta_{(i\bullet)}$, then $X_i$ is a general member of a nef divisor $|\sum_{j=l+1}^{r_i}D_{\rho_{j}}|$. We may obtain Givental's LG model as in (\ref{eq:givental mirror hybrid}). However, in general, the collection $\{\rho_1\}, \dots, \{\rho_l\}, \{\rho_j|j=l+1, \dots, r_i\}$ does not form a nef partition of the polar dual $\Delta^{(i)}$ (here we abusely write $\rho_\bullet$ for vertices of $\Delta^{(i)}$). Even if $l=1$ and the semi-stabe degeneration is Tyurin degeneration, we still have the same issue. Not only that, the mirror LG model for each component are not necessarily algebraic. 
        
    \item Another issue is that the procedure of taking the polydisks in $\bbP^l$ and completing the fibration $h_i:Y_i \to \Delta_i$ is not algebraic. Therefore, it is difficult to compute invariants like Hodge numbers to do the consistency check for being mirror to the pair $(X_i, \cup_{j \neq i}X_{ij})$. 
     \end{enumerate}

\section{Poincaré duality of hybrid Landau-Ginzburg models}\label{sec:Appendix}

Let $(Y, \omega, h=(h_1, \dots, h_N):Y \to \C^N)$ be a hybrid LG model of rank $N$. As the Kähler form $\omega$ does not play a crucial role in the following discussion, we drop it from the notation. We recall the definition of the compactified hybrid LG model of $(Y, h:Y \to \C^N)$ \cite[Section 5]{mypapermirror}.

\begin{defin}\label{def:compactified LG general}
	A compactified hybrid LG model of $(Y,h:Y \to \C^N)$ is a datum $((Z,D_Z), f:Z \to (\bbP^1)^N)$ where:
		\begin{enumerate}
			\item $Z$ is a smooth projective variety and $f=(f_1, \dots, f_N) :Z \to (\bbP^1)^N$ is a  projective morphism where each morphism $\Hat{f_i}=(f_1, \dots, \Hat{f_i}, \dots, f_N):Z \to (\bbP^1)^{N-1}$ is the compactification of the potential $\Hat{h_i}=(h_1, \dots, \hat{h_i}, \dots, h_N):Y \to \C^{N-1}$ for all $i=1, \dots, N$;
			\item the complement of $Y$ in $Z$ is a simple normal crossing anti-canonical divisor $D_Z:=D_{f_1} \sqcup \cdots \sqcup D_{f_N}$ where
    	$D_{f_i}:=(f_i^{-1}(\infty))_{red}$ is the reduced pole divisor of $f_i$ for all $i=1, \dots, N$;
    	\item the morphism $f:Z \to (\bbP^1)^N$ is semi-stable at $(\infty, \dots, \infty)$;
		\end{enumerate}
		In particular, we call the compactified LG model $((Z, D_Z), f)$ \textit{tame} if the pole divisor $f_i^{-1}(\infty)$ is reduced for all $i$.  
	\end{defin}
	 
Consider the logarithmic de Rham complex $(\Omega_Z^\bullet(\log D_Z), d)$. We define $f$-adapted de Rham complex of $Z$, denoted by $(\Omega_Z^\bullet(\log D_Z, f),d)$, to be subcomplex which is preserved by the wedge product of all $df_i$'s.
\begin{equation*}
    \Omega_Z^a(\log D_Z, f):=\{\eta \in \Omega_Z^a(\log D_Z)|\eta \wedge df_i \in \Omega^{a+1}_Z(\log D_Z) \text{ for all } i=1,\dots,n\}
\end{equation*}

First, note that it is a locally free $\cO_Z$-module of rank $\binom{n}{a}$. Here is the local description. Denote $D_{\infty}$ the corner that is the intersection of $D_{f_i}$'s. For $p \in D_{\infty}$, we can find local analytic coordinates $z_{i1}, \cdots, z_{ik_i}$ centered at $p$ with $k_1+ \cdots k_N \leq n$ such that in a small neighborhood of $p$, the divisor $D_{f_i}$ is given by $\prod_{i=1}^{k_1}z_i=0$ and the potential $f_i$ is given by 
	\begin{equation*}
		f_i(z_1, \cdots, z_n)=\frac{1}{z_{i1}^{a_{i1}}\cdots z_{1k_1}^{a_{ik_i}}}
	\end{equation*}
	for some $a_{ik_i} \geq 1$.
	\begin{lem} \label{lem: local description f-adapted}
	The $f_i$-adapted de Rham complex $\Omega^a_Z(\log D_Z, f_i)$ is locally free of rank $\binom{n}{a}$ for all $a \geq 0$. Explicitly,
	\begin{equation*}
	\Omega_Z^a(\log D_Z, f_i)=\bigoplus_{p=0}^a\left[\frac{1}{f_i} \wedge^pW_i \oplus d\log f_i \wedge \left(\wedge^{p-1}W_i\right)\right]\bigotimes \wedge^{a-p}R_i
	\end{equation*}
	where $W_i$ is spanned by logarithmic 1-forms associated to the vertical part of $f_i:Z \to \bbP^1$ and $R_i$ is spanned by holomorphic 1-forms on $Y$ and logarithmic 1-forms associated to the horizontal part of $f_i:Z \to \bbP^1$.
	\end{lem} 
\begin{proof}
    See \cite[Lemma 2.12]{KKPbogomolov}
\end{proof}
	The above local description allows one to describe $\Omega^a_Z(\log D_Z,f)$ for all $a \geq 0$. Explicitly, we have
	\begin{equation}\label{eq:local description}
	\begin{aligned}
	\Omega_Z^a(\log D_Z, f)=\bigoplus_{p_1+\cdots+p_N=0}^a \bigotimes_{i=1}^N  \left[\frac{1}{f_i} \wedge^{p_i}W_i \oplus d\log f_i \wedge \left(\wedge^{p_i-1}W_1\right)\right] \wedge^{a-(p_1+\cdots+p_N)}R
	\end{aligned}
	\end{equation}
	where $R$ is spanned by holomorphic 1-forms on $Y$. Consider the cup(=wedge) product 
	\begin{equation*}
	    \cup:\Omega^a_Z(\log D_Z, f) \otimes \Omega^{n-a}_Z(\log D_Z, f) \to \Omega^n_Z(\log D_Z,f)=\Omega^n_Z(\log D_Z)
	\end{equation*}
	Note that this is non-degenerate. From the description of the local form, one can see that it factors through $\Omega_Z^n$, the sheaf of holomorphic $n$-forms on $Z$. In other words, we have the following commutative diagram 
	\begin{equation*}
	    \begin{tikzcd}
	         \Omega^a_Z(\log D_Z, f) \otimes \Omega^{n-a}_Z(\log D_Z, f) \arrow[r, "\cup"] \arrow[rd]&  \Omega^n_Z(\log D_Z) \\
	         &\Omega^n_Z \arrow[u, hook]
	    \end{tikzcd}
	\end{equation*}
	This provides the natural isomorphism of locally free $\cO_Z$-modules 
	\begin{equation*}
	    \Omega_Z^a(\log D_Z, f) \cong Hom_{\cO_Z}(\Omega_Z^{n-a}(\log D_Z, f), \Omega_Z^n) \cong (\Omega_Z^{n-a}(\log D_Z, f))^* \otimes \Omega_Z^n
	\end{equation*}
	Therefore, we have a perfect pairing which is given by the composition of cup product with the natural trace map:
	\begin{equation*}
	    \bbH^q(Z, \Omega_Z^p(\log D_Z, f)) \otimes \bbH^{n-q}(Z, \Omega_Z^{n-p}(\log D_Z, f)) \to H^n(X, \Omega_Z^n) \xrightarrow{Tr} \C 
	\end{equation*}
    hence we have the Serre duality for the $f$-adapted de Rham forms. 
\begin{prop}\cite[Proposition 5.25]{mypapermirror}
The $f$-adapted de Rham complex $\Omega^\bullet_Z(\log D_Z,f)$ satisfies the Hodge-to de Rham degeneration property. In particular, we have 
\begin{equation*}
    \bbH^a(Z, \Omega_Z^\bullet(\log D_Z, f)) \cong \bigoplus_{p+q=a}\bbH^q(Z, \Omega_Z^p(\log D_Z, f))
\end{equation*}
\end{prop}

It is known that the cohomology $\bbH^a(Z, \Omega_Z^\bullet(\log D_Z, f))$ is isomorphic to the relative cohomology $H^a(Y, \sqcup_{i=1}^N Y_i,\C)$ where $Y_i$ is a smooth generic fiber of $h_i:Y \to \C$ near the infinity and $\sqcup_{i=1}^N Y_i$ is the normal crossing union of $Y_i$'s. The gluing property of the hybrid LG model (Proposition \ref{prop:Gluing Property in general}) provides an isomorphism $H^a(Y, \sqcup_{i=1}^N Y_i,\C) \cong H^a(Y, Y_{sm}, \C)$ where $Y_{sm}$ is a smooth generic fiber of the associated ordinary LG potential $\sw:Y \to \C$. Therefore, we have

\begin{thm}\label{thm:PD hybrid}(Poincaré duality)
    Let $(Y, h:Y \to \C^N)$ be a rank $N$ hybrid LG model. Then for $a \geq 0$, there is an isomorphism of cohomology groups 
    \begin{equation*}
        H^a(Y, Y_{sm},\C) \cong H^{2n-a}(Y, Y_{sm},\C)^*
    \end{equation*}
    where $n=\dim_\C Y$.
\end{thm}

\begin{rem}
    There is a secret sign issue. When we say Poincaré duality, we want to take the topological trace map instead of the algebraic trace map. This amounts to multiplying the sign $(-1)^n$. 
\end{rem}

\printbibliography
\end{document}